\documentclass[11pt,leqno]{amsart}

\newtheorem{prop}{Proposition}
\newtheorem{theo}{Theorem}
\newtheorem{Lemma}{Lemma}
\newtheorem{cor}{Corollary}

\newtheorem{rem}{Remark}

\newcommand{\na}{{D^g}}
\newcommand{\om}{\omega}
\newcommand{\Om}{\Omega}
\newcommand{\la}{\lambda}
\newcommand{\La}{\Lambda}

\newcommand{\cal}{\mathcal}
\newcommand{\id}{{\rm Id}|_{TM}}
\newcommand{\Lie}{\rm}

\newcommand{\llra}{\!\!\joinrel{\hbox to 30pt{\rightarrowfill}}}
\newcommand{\lllra}{\!\!\joinrel{\hbox to 50pt{\rightarrowfill}}}
\newcommand{\llllra}{\!\!\joinrel{\hbox to 60pt{\rightarrowfill}}}
\newcommand{\lllllra}{\!\!\joinrel{\hbox to 70pt{\rightarrowfill}}}
\newcommand{\llllllra}{\!\!\joinrel{\hbox to 75pt{\rightarrowfill}}}

\title[]{Self-dual Einstein Hermitian four manifolds}
\author{VESTISLAV APOSTOLOV
AND PAUL GAUDUCHON}
\thanks{The first author was supported in part by NSF grant INT-9903302}
\address{Vestislav Apostolov \\ D{\'e}partement de Math{\'e}matiques\\
UQAM\\ C.P. 8888 \\ Succ. Centre-ville \\ Montr{\'e}al (Qu{\'e}bec) \\
H3C 3P8 \\ Canada}
\email{apostolo@math.uqam.ca}

\address{Paul Gauduchon \\ CMAT\\ {\'E}cole Polytechnique \\ UMR 7640 du CNRS
\\ 91128 Palaiseau \\ France}
\email{pg@math.polytechnique.fr}

\begin{document}

\begin{abstract} We provide a local classification of
self-dual Einstein Riemannian 
four manifolds admitting a positively oriented Hermitian structure and
characterize  those 
which carry a hyperhermitian, non-hyperk{\"a}hler 
structure compatible with the negative 
orientation. We finally show that self-dual Einstein 4-manifolds obtained as
quaternionic quotients of the Wolf spaces ${\mathbb H}P^2$, ${\mathbb H}H^2$,
$SU(4)/S(U(2)U(2))$, and $SU(2,2)/S(U(2)U(2))$ are always Hermitian.

\vspace{0.1cm}
\end{abstract}
\maketitle
\section*{Introduction}

The main goal of this paper is to give a local description of 
all self-dual Einstein
$4$-manifolds $(M,g)$ which admit a positive Hermitian structure.

It follows from a (weak) Riemannian version of the Goldberg-Sachs theorem 
\cite{PB,Bo4,Nu,AG} that a Riemannian  Einstein $4$-manifold locally admits a
positive Hermitian structure if and only if the {\it self-dual Weyl
tensor}  $W^+$ is {\it
  degenerate}. 
This means that  
at any point of $M$ at least two of the three  eigenvalues of
$W^+$ coincide, when $W ^+$ is viewed as a symmetric traceless
operator  acting on 
the three-dimensional space of self-dual 2-forms.

Riemannian Einstein 4-manifolds with
degenerate self-dual Weyl tensor have been much studied by
A. Derdzi\'nski; we here recall the following facts taken from \cite{De}: 
\begin{enumerate} 

\item[{\rm (i)}]  $W^+$ either vanishes identically or else has no
  zero, i.e. has exactly two distinct eigenvalues at any point (one of them, say $\lambda$,
 is simple; the
 other one is
 of multiplicity $2$, and therefore equals  $- \frac{\lambda}{2}$ as 
$W ^+$ is trace-free).
\smallskip

\item[{\rm (ii)}]  In the latter case,  the K{\"a}hler form of the 
Hermitian structure $J$ is a generator of the simple eigenspace 
of $W^+$ --- in particular, the conjugacy class of $J$ is uniquely
defined by the metric ---
and the conformal
metric ${\bar g} = |W^+|^{\frac{2}{3}}g$ is   {\it K{\"a}hler} 
with respect to $J$. 
\smallskip

\item[{\rm (iii)}]  If, moreover, $g$ is  assumed to be {\it
    self-dual} --- meaning that the {\it anti-self-dual} Weyl tensor, 
  $W ^-$,  vanishes identically --- 
the simple eigenvalue $\lambda$ of $W ^+$ is constant (equivalently, the norm $|W ^+|$
is constant) if and only if $(M, g)$ is
locally symmetric, i.e., a real
or complex space form. 

\end{enumerate}

We then have a natural bijection between the following three classes of
Riemannian $4$-manifolds (see Lemma \ref{de-ga} below): 

\begin{enumerate}

\item Self-dual Einstein $4$-manifolds with degenerate self-dual Weyl
  tensor $W ^+$, such that $|W ^+|$ is not constant.

\smallskip

\item Self-dual Einstein Hermitian $4$-manifolds which are neither
conformally-flat  nor K{\"a}hler.

\smallskip

\item Self-dual K{\"a}hler manifolds with nowhere vanishing and
non-constant scalar curvature.

\end{enumerate}

In this correspondence, the Riemannian metrics are defined on the  
same manifold and
belong to the  same conformal class. 
Observe that each class is defined by an algebraic closed condition (the
vanishing of some tensors) and an open  genericity condition.

Since the compact case is completely understood,  see e.g. 
\cite{BYC} or \cite{De,Bu,It,Bo4, ADM} for a classification,   
the paper will concentrate on the
local situation.  

The first known examples of (non-locally-symmetric) 
self-dual Einstein Hermitian metrics have been  metrics of
cohomogeneity one 
under the isometric action of a four-dimensional Lie
group. Einstein metrics  which are of cohomogeneity one under the action of a 
four-dimensional Lie group are automatically Hermitian \cite{De}. By
using this remark, A. Derdzi\'nski   constructed \cite{de} 
a family of  cohomogeneity-one
self-dual Einstein Hermitian metrics under the action of 
${\mathbb R}\times {\rm Isom}({\Bbb R}^2)$, U(1,1) and U(2);
this family actually includes (in a rather implicit way) 
the well-known {\it Pedersen-LeBrun metrics}
\cite{P,Le} which play an important r{o}le in 
Section 3 of this paper.

It is {\it a priori} far from obvious  that there are any 
other examples of self-dual Einstein Hermitian 4-manifolds, since the conditions of being self-dual, Einstein and
Hermitian constitute an over-determined second order PDE system for the metric 
$g$. We show however that there 
are actually many other examples; more precisely, we 
classify all local solutions of this system 
and provide  a simple,  {\it explicit} (local) Ansatz for self-dual Einstein
Hermitian 4-manifolds  (see Theorem \ref{th3} and Lemma \ref{integrate} for a precise
statement).

An amazing, a priori  unexpected  fact comes out from the argument and
explains a posteriori the integrability of the above mentioned
Frobenius system~: {\it 
all self-dual Einstein Hermitian 
metrics admit  
a local isometric action 
of ${\Bbb R}^2$ with two-dimensional orbits}  (Theorem \ref{th3} and Remark 3).
In particular, these metrics locally  fall into  the more general
context  of
self-dual metrics with torus action considered in \cite{joyce}
and, more recently, in \cite{Ca0,Ca1} (see Remark \ref{rem3} (ii)).

It turns out that this property of having more
(local) 
symmetries than expected is actually shared by K{\"a}hler metrics with
vanishing Bochner tensor in all
dimensions, as shown in the recent work  of R. Bryant
\cite{Br} (see \cite{Br} for precise statements). 
Since the Bochner tensor of a K{\"a}hler
manifold of real dimension four is the same as the anti-self-dual tensor $W ^-$ --- so that
Bochner-flat K{\"a}hler metrics are a natural generalization of
self-dual K{\"a}hler metrics in higher dimensions --- by using the
correspondence given by Lemma 2, Bryant's work  provides an
alternative approach to  our classification in Section 2.

Moreover, Bryant's work includes a large section devoted to complete
metrics; in particular, by specifying his general techniques to dimension four,
he has been able (again via Lemma 2)  to give {\it complete} examples
of self-dual Einstein
Hermitian 4-manifolds, corresponding to the generic case considered
in Theorem \ref{th3}.

\smallskip

The paper is organized as follows:  

\smallskip

Section 1 displays  the background material; the notation closely
follows our previous work \cite{AG} --- with the exception of the Lee
form, whose definition here is slightly different ---  and we send
back  the reader to
\cite{AG} for more details and references.

\smallskip

Section 2.1 provides a complete description  of  (locally defined)
cohomogeneity-one
self-dual Einstein Hermitian metrics (Theorem \ref{th2}). 
It turns out that they all
admit
a local isometric action (with three-dimensional orbits) 
of certain four-dimensional Lie groups,
such that  the metrics
can be put in a {\it diagonal} form; in other
words,  they are
{\it biaxial diagonal Bianchi metrics of type A}, see {e.g.} \cite{Tod1,CP}.
Theorem \ref{th2} relies on  the  fact 
that 
every (non-locally-symmetric) 
self-dual Einstein Hermitian metric $(g,J)$ 
has a distinguished non-trivial Killing field,  namely 
$K= J{\rm grad}_g(|W^+|^{-\frac{1}{3}})$, \cite{De}. Then, the
Jones-Tod  reduction  with respect to $K$ \cite{Tod2} provides 
a three-dimensional space of {\it constant curvature}. 
The diagonal form of the metrics follows from \cite{Tod2} and
\cite{Tod1} (a unified presentation for these  
cohomogeneity-one metrics also appears in \cite{CP}).
To the best of our knowledge, apart from these metrics no other examples of 
self-dual Einstein Hermitian metrics were known in the literature (see
however Section 4).
\smallskip

Section 2.2  is devoted to the generic case,  when the
metric is neither locally-symmetric nor of cohomogeneity one. 
Our approach is similar to  Armstrong's one  in \cite{Arm}:
When considering the Einstein condition alone,  
the Riemannian Goldberg-Sachs theorem together
with Derdzi\'nski's results reported above  imply a number of
relations for 
the 4-jet of an Einstein Hermitian metric
(Sec. 2.1, Proposition \ref{prop2}); these happen to be  the
only obstructions for prolonging the 3-jet solutions of the problem to
4-jet and no further obstructions appear when  
reducing the equations for non-K{\"a}hler, 
non-anti-self-dual Hermitian 
Einstein
4-manifolds to a  (simple) perturbated 
${\rm SU}(\infty)$-Toda field equation \cite{Arm, Pl-Pr}. 
If, moreover, we insist  that $g$ be   also {\it self-dual}, we find 
further relations for  the 5-jet of the metric and we show 
that  they  have  the form of an integrable  closed
Frobenius system of PDE's for 
the parameter space of  the 4-jet  of the metric. We thus prove the
local existence of non-locally symmetric and non-cohomogeneity-one 
self-dual Einstein Hermitian 
metrics (Theorem \ref{th3}).  It turns out that this 
Frobenius system can be explicitly integrated 
(Lemma 3). We thus obtain  a uniform local description 
for {\it all}   
self-dual Einstein Hermitian 
metrics in an explicit way.

Section 3  is devoted to the subclass of self-dual Einstein Hermitian 
metrics which admit a compatible, non-closed,  anti-self-dual hypercomplex
structure. This is the same, locally, as the class of self-dual Einstein Hermitian 
 metrics which admit a non-closed Einstein-Weyl connection (see
 Section 1.2). From this viewpoint, it is a particular case of
 four-dimensional conformal metrics  which admit two distinct
 Einstein-Weyl connections. In our case, one of them is the
 Levi-Civita connection of the Einstein metric, whereas the  other one
 is {\it non-closed},  hence, because of Proposition \ref{prop4}, 
 attached to a non-closed hyperhermitian structure. 
(Recall that  a conformal 4-manifold admitting two distinct
{\it closed} Einstein-Weyl structures is necessarily conformally flat
(folklore), and that, conversely, 
every conformally flat 4-manifold only admits closed Einstein-Weyl 
structures \cite{ET}, see also Proposition \ref{prop4} and
Corollary \ref{cor1} below).

It turns out that self-dual Einstein Hermitian 
 metrics which admit a compatible, non-closed,  anti-self-dual hypercomplex
structure, actually admit a second one and thus fall 
in the {\it bi-hypercomplex} situation described by Madsen in
 \cite{Mad}; in particular, these metrics admit a local action of 
${\rm U(2)}$, with three-dimensional orbits, and are diagonal Bianchi XI
 metrics, see Theorem \ref{th1} below. 

Notice that a general description of
(anti-self-dual)  metrics  admitting  two distinct compatible hypercomplex
structures appears  in \cite{Ca2}, see also \cite{BCM}, whereas a
family of self-dual Einstein metrics with compatible non-closed hyperhermitian
structures, parameterized by  
holomorphic functions  of one variable,  
has been constructed in \cite{CT}.

In Section 4, we show that all anti-self-dual,
  Einstein four dimensional {\it orbifolds}  obtained by quaternionic
  K{\"a}hler reduction from the eight dimensional quaternionic K{\"a}hler
  Wolf spaces ${\mathbb H}{P}^2$, $SU(4)/S(U(2)U(2))$
  and their non-compact duals
  (see \cite{galicki1,galicki2} and \cite{G-L}) are
  actually Hermitian with respect to the opposite orientation,
  hence locally isomorphic to  metrics described in Section 2.
  These orbifolds include the {\it  weighted projective planes}
  ${\mathbb C} P^{[p_1,p_2,p_3]}$ for
 integers $0<p_1\le p_2\le p_3$ satisfying $p_3<p_1+p_2$,
 cf. \cite[Sec. 4]{G-L}. On these  orbifolds,
 Bryant has constructed Bochner-flat
 K{\"a}hler metrics with everywhere positive scalar curvature, 
 hence also self-dual, Einstein Hermitian metrics
 according to Lemma 2 below, \cite[Sec. 4.3]{Br};
 in view of the results of Section 2,
 Galicki-Lawson's and Bryant's metrics agree locally,
 but the issue  as to whether they agree globally remains unclear.

\vspace{0.2cm}
\noindent
{\bf Acknowledgments.} The first-named author thanks the 
Dipartimento di Matematica, Universit{\`a} di Rome Tre and the 
Max-Planck-Institut in Bonn
for hospitality during the preparation
of this paper.  He would like to express his gratitude to J. Armstrong 
for explaining his approach to Einstein Hermitian metrics
and many illuminating discussions.
The authors warmly 
thank  S. Salamon for being an initiator of this work and for gently
sharing his expertise, and 
C. LeBrun, C. Boyer, K. Galicki, whose comments
are at the origin  of the last section of the paper.
It is also a pleasure for us to thank  N. Hitchin, S. Marchiafava,
H. Pedersen, P. Piccinni, M. Pontecorvo  and K.P. Tod for their interest
and  stimulating conversations, and 
R. Bryant for his interest and remarks.

Finally, a special   aknowledgment  is due to D. Calderbank for his friendly assistance
in  carefully reading the manuscript, checking computations, correcting mistakes and suggesting 
improvements; he in particular decisively contributed to  improving
the paper by pointing out a mistake in a former
version and  thus revealing the rational
character of the metrics described  in Section 2.2.

\section{Einstein metrics, Hermitian structures and Einstein-Weyl geometry
in dimension 4}

\subsection{Einstein metrics and compatible Hermitian structures}

In the whole paper $(M, g)$ denotes an oriented Riemannian
four-dimensional manifold.

A specific  feature of the four-dimensional Riemannian geometry is 
the splitting 
\begin{equation} \label{split1} AM = A ^+M \oplus A ^-M, \end{equation}
of the Lie algebra bundle, $AM$,  of skew-symmetric endomorphisms of the
tangent bundle, $TM$, into the direct sum of two Lie algebra subbundles, $A
^{\pm}M$, derived from the Lie algebra splitting $\mathfrak{so}
(4) =  \mathfrak{so}
(3) \oplus \mathfrak{so}(3)$ of the orthogonal Lie algebra  $\mathfrak{so}
(4)$ into the direct sum of two copies of $\mathfrak{so}
(3)$. 

A  similar decomposition occurs for the bundle  $\Lambda ^2M$  of
$2$-forms
\begin{equation} \label{split2} \Lambda ^2 M = \Lambda  ^+M \oplus
  \Lambda  ^-M, \end{equation} 
given by the spectral decomposition of
the Hodge-star operator, $*$, whose restriction to $\Lambda ^2 M$ is an
involution; here,  
$\Lambda  ^{\pm} M$ is  the eigen-subbundle for the eigenvalue
$\pm$ of $*$. 

Both decompositions are actually determined by the
conformal metric $[g]$ only. When $g$ is fixed, 
$\Lambda ^2 M$ is identified to $AM$ by
setting: $\psi (X, Y) = g(\Psi (X), Y)$, for any $\Psi$ in $AM$ and any
vector fields $X, Y$; then, we can arrange signs in (\ref{split1}) so
that (\ref{split1}) and
(\ref{split2}) are identified to each other. A similar 
decomposition and a similar identification
occur for the bundle $\Lambda ^2 (TM)$ of bivectors.

Sections of $\Lambda  ^+M$, resp. $\Lambda  ^-M$, are called {\it
  self-dual}, resp. {\it anti-self-dual}, and similarly for
  sections of $AM$ or $\Lambda ^2 (TM)$.

In the sequel, the   vector bundles $AM$, $\Lambda ^2 M$ and
$\Lambda ^2 (TM)$ will be freely identified to
each other; similarly, the cotangent bundle $T ^*M$ will be freely
identified to $TM$; when no confusion can arise, the inner product determined by $g$ will be
simply 
denoted by $(\cdot, \cdot)$; we adopt the convention that  $(\Psi _1 , \Psi
_2) = - \frac{1}{2} {\rm tr} \, (\Psi _1 \circ \Psi _2)$, for sections
of $AM$, and the corresponding convention for $\Lambda ^2 M$ and $\Lambda ^2 (TM)$.

The Riemannian curvature, $R$, is defined by
$R _{X, Y} = D ^g _{[X, Y]} - [D ^g _X, D ^g _Y],$
where $D ^g$ denotes the Levi-Civita connection of $g$; $R$ is thus  a
$AM$-values $2$-form, but will be rather considered as a section of
the bundle $S^2 (\Lambda ^2M)$ of symmetric 
endomorphisms of $\Lambda ^2 M$.  

The Weyl tensor, $W$, commutes with $*$ and, accordingly,  splits as $W = W ^+ + W ^-$, where $W ^{\pm} = \frac{1}{2} (W \pm W \circ *)$;  $W ^+$
is called the {\it self-dual Weyl tensor}; it acts trivially on
$\Lambda ^- M$ and will be considered in the sequel as a
field of 
(symmetric, trace-free) endomorphisms of $\Lambda ^+ M$; similarly,
the {\it anti-self-dual Weyl tensor} $W ^-$ will be considered as a
field of endomorphismes of $\Lambda ^ - M$.

The Ricci tensor, ${\rm Ric}$, is the symmetric bilinear
form defined by ${\rm Ric} (X, Y) = {\rm tr} \, \{ Z \to R _{X, Z} Y \}$;
alternatively, ${\rm Ric} (X, Y) =
\sum _{i = 1}^4 (R _{X, e _i} Y, e _i)$ for any $g$-orthonormal
basis $\{ e _i \}$. We then have ${\rm Ric} = \frac{s}{4} \, g + {\rm Ric} _0$,
where $s$ is the scalar curvature (= the trace of ${\rm Ric}$ with respect
to $g$)  and ${\rm Ric} _0$ is the {\it trace-free Ricci tensor}. The latter can be
made into a section of $S ^2 (\Lambda ^2 M)$, then denoted by
$\widetilde{{\rm Ric} _0}$, by putting
$ \widetilde{{\rm Ric} _0} (X \wedge Y) = {\rm Ric} _0 (X) \wedge Y + X
  \wedge {\rm Ric} _0 (Y). $

It is   readily checked  that  $\widetilde{{\rm Ric} _0}$ satisfies the first Bianchi
identity, i.e. $\widetilde{{\rm Ric} _0}$ is a tensor of the same kind as $R$ itself,  as well as $W
^+$ and $W ^-$; moreover, $ \widetilde{{\rm Ric} _0}$  anti-commutes with $*$,
so that it can be viewed as a field of homomorphisms from $\Lambda ^+
M$ into $\Lambda ^- M$, or from $\Lambda ^-
M$ into $\Lambda ^+ M$ (adjoint to each other); 
we eventually get  the well-known Singer-Thorpe decomposition  of $R$, see
e.g. \cite{besse}:
\begin{equation}\label{SO(4)} R = \frac{s}{12} \, {\rm Id} _{| \Lambda ^2 M}  +
  \frac{1}{2} \widetilde{{\rm Ric}} _0 +
  W ^+ + W ^-, \end{equation}
or, in a more pictorial way
\begin{equation*} R = \begin{pmatrix} & W ^+ + \frac{s}{12} \, {\rm
  Id} _{| \Lambda ^+ M} & \frac{1}{2} \widetilde{{\rm Ric} _0} _{| \Lambda ^- M}
  \\ \\ & \frac{1}{2} \widetilde{{\rm Ric} _0} _{| \Lambda ^+ M} & W ^- + \frac{s}{12} \, {\rm
  Id} _{| \Lambda ^- M} \end{pmatrix} \end{equation*}

The metric $g$ is {\it Einstein}  if ${\rm Ric} _0 = 0$ (equivalently,
$g$ is Einstein if $R$
commutes with $*$).

The metric $g$ (or rather  the conformal class $[g]$) is  {\it
  self-dual}  if $W ^- = 0$; {\it anti-self-dual} if $W ^+ = 0$.

An {\it almost-complex structure}  $J$ is a field of automorphisms of $TM$ of
square $ - \id$. An   {\it integrable}  almost-complex structure is simply called a
{\it complex structure}.

In this paper,
the metric $g$, or its conformal class $[g]$, is fixed and we only  consider  
$g$-orthogonal almost-complex structures, i.e. almost-complex
structure $J$   satisfying the identity $g (JX, JY) = g (X, Y)$, so
that the pair $(g, J)$ is an {\it almost-Hermitian structure}; then, the associated bilinear form,
$F$, defined by $F (X, Y) = g(JX, Y)$ is a $2$-form, called the {\it
  K{\"a}hler form}. 

The pair $(g, J)$ is {\it Hermitian} if $J$ is integrable; {\it
  K{\"a}hler} if $J$ is parallel with respect to the Levi-Civita
  connection $D ^g$; if $(g, J)$ is K{\"a}hler  then $J$ is integrable and $F$
  is closed; conversely, these two conditions together imply that $(g, J)$ is
  K{\"a}hler.

A $g$-compatible almost-complex structure $J$ is either a section of
$A ^+M$ or a section of $A ^-M$; it is  called {\it positive}, or
{\it self-dual},  in
the former case, {\it negative}, or {\it anti-self-dual}  in the latter case. Alternatively, the
K{\"a}hler form {\rm is} either self-dual or  
anti-self-dual.
Conversely, any
section $\Psi$ of $A ^+M$, resp. $A ^- M$,  such that $|\Psi| ^2 = 2$,
is a positive, resp. negative, $g$-orthogonal almost-complex
structure. It follows that any non-vanishing section, $\Psi$,  of $A
^+M$ --- if any --- determines a (positive) almost-complex structure
$J$,  defined by $J = \sqrt{2} \frac{\Psi}{|\Psi|}$ (similarly for
non-vanishing sections of $A ^- M$). 

Whereas the existence of a (positive) $g$-orthogonal  almost-complex 
structure is a purely topological problem, the similar issue   for
{\it complex} structures heavily depends on the
geometry of $g$, and this dependence is essentially measured by the self-dual Weyl 
tensor  $W ^+$. 

This assertion  can be made more precise in the following way. We denote by $\lambda _ +
\geq  \lambda _0 \geq  \lambda _-$ the eigenvalues of $W ^+$ at some
point, $x$, of $M$, and we assume that $W ^+$ does not vanish at $x$; 
equivalently, since $W ^+$ is trace-free, we assume that  $\lambda _+ -
\lambda _-$ is positive; we denote by $F
_{+}$ an eigenform of $W ^+$ with respect to $\lambda _+$,
normalized by $|F _+|  ^2 = 2$; similarly, $F _-$ denotes  an eigenform
of $W ^+$ for $\lambda _-$, again  normalized by $|F _-|  ^2 = 2$; the
{\it roots}, $P$,   of $W ^+$ at $x$ are then  defined by 
$P = \frac{(\lambda_+ - \lambda_0 )^{\frac{1}{2}}}{(\lambda_+ -
  \lambda_- )^{\frac{1}{2}}}F_-  +  
\frac{(\lambda_0 - \lambda_- )^{\frac{1}{2}}}{(\lambda_+ -
\lambda_-)^{\frac{1}{2}}}F_+;$
it is easily checked that this  expression actually determine {\it
  two}
distinct pairs of opposite roots in the generic case,  when the
eigenvalues are all distinct, and {\it  one}  pair in the degenerate case, 
when  $\lambda _0$ is
equal to either $\lambda _+$ or $\lambda _-$.  

It is   a basic fact  that when $J$ is a positive, $g$-orthogonal
{\it complex}  structure defined on $M$, the value of $J$ at any point $x$
where $W ^+$ does not vanish must be equal to a root of $W ^+$  at that
point. This means that on the open subset of $M$ where $W ^+$ does not
vanish, the conjugacy class of a positive, $g$-orthogonal
complex structure --- if any --- is almost entirely determined by $g$
(in fact by $[g]$), with at most a $2$-fold ambiguity. 

On the other hand, it is an easy
consequence of the integrability theorem in  \cite{AHS} that $A ^+M$
can be locally trivialized by integrable (positive, $g$-orthogonal)
almost-complex structures if and only if $[g]$ is anti-self-dual.

In the sequel, $W ^+$ will be called {\it degenerate} at some point $x$ if it has at
most two distinct eigenvalues at that point. The terms {\it
  anti-self-dual} and {\it non-anti-self-dual} will be abbreviated as
ASD and non-ASD respectively.

For a given non-ASD metric $g$ it is a subtle question to decide whether
the roots of  $W ^+$ actually provide complex structures (this is of
course not true in general). The situation is quite different  if $g$
is Einstein. It is then settled by the following  (weak) Riemannian version
of the Goldberg-Sachs theorem, cf. \cite{De,PB,Nu,AG}:
\begin{prop}\label{prop1} Let $(M,g)$ be an oriented Einstein
4-manifold; then the following three conditions are equivalent:
\begin{enumerate}
\item[(i)] $W^+$ is everywhere degenerate;
\smallskip

\item[(ii)] there exists a positive $g$-orthogonal complex structure in a neighbourhood
of each point of $M$;
\smallskip

\item[(iii)] $(M,g)$ is either ASD or   $W^+$ has 
two distinct eigenvalues at each  point.
\end{enumerate}
\end{prop}

A consequence of this proposition  is  that the self-dual Weyl tensor $W
^+$ of a non-ASD Einstein Hermitian $4$-manifold 
nowhere vanishes and has  two distinct eigenvalues at any point, one
simple, the other one of multiplicity $2$; moreover, the K{\"a}hler form $F$ is an 
eigenform of $W ^+$ for the simple eigenvalue. Conversely, 
for any oriented, Einstein $4$-manifold 
whose $W ^+$  has  two distinct eigenvalues,  
the generator of the simple eigenspace of
$W^+$ determines a (positive)  Hermitian structure.

For any positive $g$-orthogonal almost-complex structure $J$, $A ^+ M$
splits as follows:
\begin{equation} \label{splitJ1} A ^+ M = {\mathbb R}\cdot {J} \oplus A ^{+, 0} M,
\end{equation}
where ${\mathbb R}\cdot{J}$ is the trivial subbundle generated by
$J$ and $A ^{+, 0} M$ is the orthogonal complement (equivalently, $A
^{+, 0} M$ is the subbundle of elements of $A ^+ M$ that anticommute with $J$); $A
^{+, 0} M$ is a rank $2$ vector bundle and will be also considered as
a complex line bundle by putting $J \Phi = J \circ \Phi$. We have the
corresponding decomposition  
\begin{equation} \label{splitJ2} \Lambda  ^+ M = {\mathbb R}\cdot{F} \oplus \Lambda
  ^{+, 0} M,
\end{equation}
where $\Lambda
  ^{+, 0} M$ is the subbundle of $J$-anti-invariant $2$-forms,
  i.e. $2$-forms satisfying $\phi (JX, JY) = - \phi (X, Y)$; again, $\Lambda
  ^{+, 0} M$ is  considered as a complex line bundle by putting $(J
  \phi) (X, Y) = - \phi (JX, Y) = - \phi (X, JY)$. As complex line
  bundles, both $A ^{+, 0} M$ and $\Lambda
  ^{+, 0} M$ are identified to the {\it anti-canonical bundle} $K
  ^{-1} M = \Lambda ^{0, 2} M$ of the (almost-complex) manifold $(M,
  J)$. 

For an Einstein, Hermitian $4$-manifold, the action of $W ^+$
preserves the decompositions
(\ref{splitJ1}) and (\ref{splitJ2}).

\bigskip

The {\it Lee form} of an almost-Hermitian structure
$(g, J)$ is the real $1$-form, $\theta$, defined by
\begin{equation} \label{lee} {\rm d} F = - 2 \theta \wedge F; \end{equation}
equivalently, $\theta = - \frac{1}{2} J \,  \delta F$, where $\delta$
denotes the co-differential with respect to $g$ (here, and henceforth,
the action of $J$ on $1$-forms is defined  via the identification $T
^*M \simeq  TM$ given by the metric; we thus have $(J \alpha ) (X) =
- \alpha (JX)$, for any $1$-form $\alpha$). The reason 
for the choice of the factor
$-2$ in (\ref{lee}) will be clear  in the next subsection (notice that  a
different 
normalization is used  in our previous work \cite{AG}). 

When  $(g, J)$ is Hermitian, it is K{\"a}hler if and only if $\theta$
vanishes identically; it is conformally K{\"a}hler if  and only
if $\theta$ is exact, i.e. $\theta = - 
{\rm d}  \ln{f}$ for a positive smooth real function $f$ (then, $J$
is K{\"a}hler with respect to the
conformal metric $g' = f ^{-2} \, g)$; it is locally conformally
K{\"a}hler  ---
lcK for short --- if and only if $\theta$ is closed, hence locally of
the above type. 

The Lee form clearly satisfies $({\rm d} \theta, F) = 0$; this means
that the self-dual part, ${\rm d}
  \theta  ^+$,  of ${\rm d}
  \theta $ is a section of the rank $2$ subbundle,
  $\Lambda ^{+,0} M$.

In the Hermitian case, ${\rm d}
  \theta  ^+$ is an eigenform of $W ^+$ for the mid-eigenvalue $\lambda
  _0$; moreover, $\lambda _0 = - \frac{\kappa}{12}$, where $\kappa$ is the {\it
  conformal scalar curvature}, of which a more direct definition is
  given in the next subsection; ${\kappa}$ is related to the (Riemannian) scalar curvature $s$ by 
\begin{equation} \kappa  = s +
  6 \, (\delta \theta - |\theta| ^2),  \end{equation}
and we also have
\begin{equation} \kappa = 3 \, (W ^+ (F), F), \end{equation} 
see \cite{Va2,Ga1}. Notice that, in the Hermitian case, the mid-eigenvalue
$\lambda _0$ of $W ^+$ is always a {\it smooth} function (this,
however,  is not
true in general for the remaining  two eigenvalues of $W ^+$,
$\lambda _+$ and $ \lambda _-$, which   are  given by:
$$\la_{\pm}=\frac{1}{24}\kappa \pm \frac{1}{8}(\kappa^2 +
32|{\rm d} \theta ^+|^2)^{\frac{1}{2}},$$
cf. \cite{AG}).

It follows that for Hermitian 4-manifolds 
the following 
three conditions are equivalent (cf. \cite{Bo2,AG}):
\smallskip

\begin{enumerate}
\item[(i)] ${\rm d} \theta ^+=0$;
\smallskip

\item[(ii)] $W^+$ is degenerate;
\smallskip

\item[(iii)] $F$ is an eigenform of $W^+$. 
\end{enumerate}

\noindent
(In the latter case $F$ is actually an eigenform for
the simple eigenvalue of $W ^+$, which  is then equal to
$\frac{\kappa}{6}$, {\rm also} equal to $\lambda _+$ or $\lambda _-$
according as  $\kappa$ is
positive or negative).
If, moreover,  $M$ is compact,  any one of the above three conditions is equivalent to
$(g,J)$ being locally conformally K{\"a}hler; if, in addition, the first
Betti number of $M$ is even, $(g, J)$ is then globally conformally K{\"a}hler \cite{Va1}.

By Proposition \ref{prop1} we conclude 
that for every Einstein Hermitian $4$-manifold, we have
${\rm d} \theta ^+=0$, i.e. ${\rm d} \theta$ is self-dual. 
In fact, a
stronger statement is true, see \cite[Prop.1]{AG} and
\cite[Prop.4]{De}:
\begin{prop}\label{prop2}
Let $(M,g,J)$ be an Einstein, non-ASD  Hermitian 4-manifold. 
Then the conformal scalar curvature $\kappa$ nowhere vanishes and the  Lee form $\theta$
is given by \ $\theta = \frac{1}{3}{\rm d}\ln{|\kappa|}$ {\rm (}in particular, 
$(g,J)$ is conformally K{\"a}hler{\rm )}.

If, moreover, $\kappa$ 
is not  constant, i.e. if $(g,J)$ is not K{\"a}hler,  then  
$K=J\rm{grad}_g(\kappa^{-\frac{1}{3}})$ is a non-trivial Killing
vector field 
with respect to  $g$, holomorphic
with respect to $J$.
\end{prop}

\subsection{Einstein-Weyl structures and anti-self-dual conformal
  metrics}

Another specific feature of the four-dimensional geometry is that  to
each conformal Hermitian structure $([g], J)$ is canonically attached a unique
{\it Weyl connection} $D$ such that $J$ is parallel with respect to
$D$; in other words, any Hermitian structure is ``K{\"a}hler'' in the
extended context of 
Weyl structures (of course, $(g, J)$ is K{\"a}hler in the usuel sense
--- the only one used in this paper --- if and only if $D$ is the
Levi-Civita connection of some metric in the conformal class $[g]$).

Recall that, given a conformal metric $[g]$,  a Weyl connection (with
respect to $[g]$) is a torsion-free linear connection, $D$,  on $M$
which preserves $[g]$; the latter condition can be reformulated as
follows: for any metric $g$ in $[g]$, there exists a real $1$-form
$\theta _g$ such that $D g = - 2 \theta _g \otimes g$; $\theta _g$ is
called the {\it Lee form} of $D$ with respect to $g$; then, the Weyl
connection $D$ and the Levi-Civita connection $D ^g$ are related by $D
= D ^g + \tilde{\theta} _g$, meaning
\begin{equation}\label{D^J} D _X Y = D ^g _X Y + \theta _g (X) Y + \theta _g (Y)
  X - g(X, Y) \, \theta _g ^{\sharp _g}, \end{equation}
where $\theta _g ^{\sharp _g}$ is the Riemannian dual of $\theta _g$ with
respect to $g$.
If $g' = f ^{-2} g$
is another metric in $[g]$, the Lee form, $\theta
_{g'}$, of $D$ with respect to $g'$ is related to $\theta _g$ by
$\theta _{g'} = \theta _g  + {\rm d} \ln{f}$.

A Weyl connection $D$ is the Levi-Civita connection of some  metric in
the conformal class $[g]$ if and only if its  Lee form with respect to
any metric $g$ in $[g]$ is exact, i.e. $\theta _g = - {\rm
  d} \ln{f}$; then, $D = D ^{f ^{-2} g}$; such a Weyl connection is
called {\it exact}. More generally, a Weyl connection is said to be
{\it closed} if its Lee form with respect to any metric in $[g]$ is
closed; then, $D$ is locally of the above type, i.e. locally the
Levi-Civita connection of a (local) metric in $[g]$. 

The definitions of the  curvature $R ^D$ and the Ricci tensor ${\rm Ric} ^D$ 
of a Weyl
connection $D$ are formally identical as the ones we gave for $D ^g$ 
(notice that the derivation
of ${\rm Ric} ^D$ from $R ^D$ requires no metric); however, $R ^D$ is
now a $AM \oplus {\mathbb R} \, {\id}$-valued $2$-form, i.e. has a
{\it scalar
part} equal to $F ^D \otimes{\id}$, where  the real $2 $-form $F
^D$, 
the so-called {\it Faraday tensor} of the Weyl connection, is equal to $- {\rm d} \theta _g$ for any
metric $g$ in $[g]$; moreover, ${\rm Ric} ^D$ is not  symmetric in general:
its  skew-symmetric part 
is equal to $\frac{1}{2} F ^D$; ${\rm Ric}^D$ is thus symmetric  if and only if $D$ is closed.

A Weyl connection $D$  is called {\it Einstein-Weyl} if the symmetric,
trace-free part of ${\rm Ric} ^D$ vanishes; with respect to any metric
$g$  
in    $[g]$, and by writing $\theta$ instead of $\theta _g$, this conditions reads
\begin{equation} \label{EW} D ^g \theta - \theta \otimes \theta + \frac{1}{4}
  (\delta \theta + |\theta| ^2) \, g - \frac{1}{2} {\rm d} \theta -
  \frac{1}{2} {\rm Ric} _0 = 0, \end{equation}
see {e.g.} \cite{Ga2}; 
for a fixed metric $g$, (\ref{EW}) should be  considered as an equation for an
  unknown $1$-form $\theta$.

The {\it conformal scalar curvature} of $D$  with respect to $g$,
denoted by $\kappa _g$, is  
the trace of ${\rm
Ric}^D$ with respect to $g$; it is related to the (Riemannian) scalar curvature $s$ by:
\begin{equation}\label{kappag}
\kappa _g = s + 6 \, (\delta \theta  - |\theta |^2),
\end{equation}
see { e.g.} \cite{Ga2}.
\smallskip

A key observation is that the Lee form, $\theta$, of an
almost-Hermitian structure $(g, J)$ is also the Lee form with respect
to $g$  of the Weyl
connection canonically attached to the conformal
almost-Hermitian structure $([g], J)$; in other words, the Weyl
connection $D$ defined by $ D = D ^g + \tilde{\theta}$ is actually
independent of $g$ in its conformal class $[g]$. The Weyl connection
$D$ defined in this way is called the {\it canonical Weyl connection} of the
(conformal) almost-Hermitian structure $([g], J)$.

The scalar curvature $\kappa _g$ of $D$ with respect to $g$ is called
the {\it conformal scalar curvature} of $(g, J)$; it coincides with
the function $\kappa$ introduced in the previous paragraph. 

The canonical Weyl connection is an  especially interesting object when
$J$ is integrable, because of the following lemma:

\begin{Lemma} \label{weyl} {\rm (i)} $J$ is integrable if and only if $DJ = 0$.

{\rm (ii)} If $J _1$ and $J _2$ are two $g$-orthogonal complex
structures, the corresponding canonical connections $D ^1$ and $D ^2$
coincide if and only if the scalar product $(J _1, J_2)$ is constant.
\end{Lemma}
\begin{proof} (i) The condition $DJ = 0$ reads 
\begin{equation} \label{integrable} D ^g  _X J  = [X \wedge \theta, J]; \end{equation}
this identity is proved e.g. in \cite{Ga1,Va2}. 

(ii) Let $p$ denote the
  {\it angle function} of $J _1$ and $J_2$, defined by $p = -
  \frac{1}{4} {\rm tr} \, (J _1 \circ J _2) = \frac{1}{2} (J _1,
  J_2)$; we then have
\begin{equation} J_1 \circ J _2 + J _2 \circ J_2 = - 2p \, {\id}. \end{equation}
Let $\theta _1$ and $\theta _2$ be the Lee forms of $D ^1$, $D ^2$;
from (\ref{integrable}) applied to $J _1$, we infer $(D ^g J _1, J_2)
=  ([J_1, J_2] X, \theta _1)$; similarly, we have $(D ^g J _2, J_1)
= ([J_2, J_1] X, \theta _2)$; putting together these two
identities, we get
\begin{equation} {\rm d} p = - \frac{1}{2} [J_1, J_2] (\theta _1 -
  \theta _2). \end{equation}
This obviously implies ${\rm d} p = 0$ if $D ^1 = D ^2$; the converse
is also true, as the commutator $[J _1, J_2]$ is invertible at each point where $J _2 \neq \pm J_1$.
\end{proof}

\smallskip

An {\it almost-hypercomplex  structure}  is the datum of three
almost-complex structures, $I _1, I_2, I_3$, such that
$$ I _1 \circ I _2 = - I _2 \circ I _1 = I _3.$$

Since $M$ is a four-dimensional manifold, any almost-hypercomplex
structure $I _1, I_2, I_3$ determines a conformal class $[g]$ with
respect to which each $I_i$ is orthogonal: $[g]$ is defined by
decreeing 
that, for any non-vanishing (local)  vector field $X$, the frame $X, I_1X,
I_2X, I_3X$ is (conformally) orthonormal; for any $g$ in the
conformal class defined in this way, we thus get an {\it
  almost-hyperhermitian structure} $(g, I_1, I_2,I_3)$; notice that
the $I_i$'s are pairwise orthogonal with respect to $g$, so that $I
_1, I_2, I_3$ is a (normalized) orthonormal frame of $A ^+ M$;
conversely, for a given Riemannian metric $g$ any  (normalized)
orthonormal frame of $A ^+ M$ is an almost-hypercomplex structure and,
together with $g$ form an almost-hyperhermitian structure.

An almost-hyperhermitian structure $(g, I_1, I_2,I_3)$ is called 
{\it hyperhermitian} if all $I_i$'s are integrable; it is called {\it
  hyperk{\"a}hlerian} if the $I_i$'s are all parallel with respect to
the Levi-Civita connection $D ^g$.

In the hyperhermitian case the canonical Weyl connections, $D ^1, D
^2, D ^3$,  of the
almost-Hermitian structures $(g, I_1)$, $(g, I_2)$, $(g, I_3)$ are the
same by 
Lemma \ref{weyl}; the common Weyl connection, $D$, is
called the {\it canonical Weyl connection}
of the hyperhermitian structure.
 
Conversely, the
condition $D ^1 = D ^2 = D ^3$ implies that $(g, I_1, I_2,I_3)$ is
hyperhermitian (this observation is due to S. Salamon and  F. Battaglia,
see  e.g.  \cite{GT}).

The canonical Weyl connection of a hyperhermitian structure $(g, I_1,
I_2,I_3)$ is
closed if and only if $I_1, I_2,I_3$ is locally hyperk{\"a}hler with
respect  to some (local) metric belonging to the conformal class $[g]$; for
brevity, a hyperhermitian structure will be called {\it closed} or
{\it non-closed} according as its  canonical Weyl connection being
closed or non-closed.

\begin{rem} {\rm In general, for any  given hypercomplex structure $I _1,
    I_2, I _3$ on a $n$-dimensional manifold, there exists a {\it unique}
    torsion--free linear connection on $M$ that preserves the $I
    _i$'s, called the {\it Obata connection}; the canonical connection thus coincides with the
    Obata connection; for $n > 4$ however, there is no conformal
    metric canonically attached to $I _1,
    I_2, I _3$ and, in general,  the Obata connection is
    not a Weyl connection.}
\end{rem}

If $(g, I_1, I_2,I_3)$ is
hyperhermitian, we have  $D I_1 = D I_2 = D I_3 = 0$, where $D$ is
the canonical Weyl connection acting on sections of $A ^+ M$; it
follows   that the connection of  $A ^+M$ induced by $D$
is {\it flat}; conversely, if $D$ is a Weyl connection, whose induced
connection on $A ^+ M$ is flat, then $A ^+ M$ can be locally trivialized by
a $D$-parallel (normalized) orthonormal frame $I_1, I_2,I_3$, which,
together with $g$, constitute a hyperhermitian structure.

The curvature, $R ^{D, A ^+M} $, of the induced connection is given  by $R ^{D, A ^+M} _{X, Y} \Psi = [R ^D _{X, Y},
\Psi]$, where $R ^D _{X, Y}$ is understood as a field of
endomorphisms of $TM$ --- more precisely  a section of $A M \oplus {\mathbb R}
\, {\id}$ ---  and $[R ^D _{X, Y},
\Psi]$ is the commutator of $R ^D _{X, Y}$ and $\Psi$; we easily
infer that the vanishing of $R ^{D, A ^+M} $ is equivalent to the
following four  conditions:

\begin{enumerate} \label{swann}

\item $W ^+ = 0$;

\smallskip

\item  $(F ^D) ^+ = 0$; if $\theta$ denotes the Lee form of $D$, this
  also reads ${\rm d} \theta ^+ = 0$;
\smallskip

\item  $D$ is Einstein-Weyl, i.e. the Lee form $\theta$ is solution of
  (\ref{EW});
\smallskip

\item The scalar curvature of $D$ vanishes identically; in view of
  (\ref{kappag}), this condition reads 
\begin{equation}\label{hypherm}
s = 6 \, (-\delta \theta  + |\theta|^2).
\end{equation}

\end{enumerate}

It follows from this discussion that,  for an ASD Riemannian
4-manifold, the existence of a compatible hypercomplex structure is
locally equivalent to the existence of an Einstein-Weyl connection
satisfying the above conditions 2 and 4 (cf. \cite{PS} or \cite{GT}). 
In this correspondence, conformally hyperk{\"a}hler  
structures  correspond to closed Einstein-Weyl structures. The existence of a
non locally hyperk{\"a}hler,  hyperhermitian structure is actually (locally)
equivalent to the existence of a non-closed Einstein-Weyl connection, 
in view of the following result of D. Calderbank:
\begin{prop}\label{prop4}{\rm (\cite{Ca})} 
Let $(M,[g],D)$ be an anti-self-dual Einstein-Weyl 4-manifold. Then either
$D$ is closed, or else $D$ satisfies conditions 2 and 4 above, i.e. is the canonical Weyl connection of a
hyperhermitian structure.
\end{prop}

Notice that in the case when $M$ is
compact, $d \theta ^+ = 0$ implies ${\rm d}\theta=0$, hence {\it any}  hyperhermitian
structure  is locally
conformally hyperk{\"a}hler; a complete classification  appears in \cite{Bo3}.
\vspace{0.2cm}

\section{Self-dual Einstein Hermitian 4-manifolds}

By Proposition \ref{prop2},  a Hermitian, Einstein 4-manifold, whose
self-dual Weyl tensor $W ^+$ has constant eigenvalues is either anti-self-dual or K{\"a}hler-Einstein, \cite{De}. 
If, moreover, the metric $g$ is self-dual,  this happens precisely when 
$g$  is locally-symmetric, i.e. when $(M, g)$ is a real or a complex space form, see
 \cite{TV}. More generally, a self-dual Einstein 4-manifold is
 locally-symmetric if and only if $W ^+$ is degenerate, with constant
 eigenvalues, \cite{De}.

In the opposite case, we have the following lemma:
\begin{Lemma} \label{de-ga} Non-locally-symmetric 
self-dual Einstein Hermitian metrics are 
in one-to-one correspondence with self-dual K{\"a}hler metrics of nowhere 
vanishing and non-constant scalar curvature.
\end{Lemma}
\begin{proof} Every self-dual
Einstein Hermitian 4-manifold $(M,g,J)$ of 
non-constant curvature is conformally related (via Proposition \ref{prop2})
to a
self-dual K{\"a}hler metric ${\bar g}$ of nowhere vanishing scalar curvature. 
A self-dual K{\"a}hler metric is locally-symmetric if and only if its scalar 
curvature is constant \cite{De}; thus, the one direction in the correspondence 
stated in the lemma follows by observing that 
${\bar g}$ is locally-symmetric as soon as $g$ is. 
Since the Bach tensor of
a self-dual metric vanishes \cite{Gau3}, 
it follows from \cite[Prop.4]{De}
that any self-dual K{\"a}hler metric of nowhere vanishing scalar curvature
gives rise to an 
Einstein Hermitian metric in the same conformal class. 
\end{proof}

In the remainder of this section, $(M, g, J)$ is an Einstein,
self-dual Hermitian $4$-manifold, and we assume that $g$ is {\it not}
locally-symmetric; in particular, $W ^+$ is degenerate, but its
eigenvalues, $\lambda, - \frac{\lambda}{2}$, or, equivalently, its
norm $|W^+| = \sqrt{\frac{3}{2}} \,  |\lambda|$, are not constant.

Since $(M,g,J)$ is not K{\"a}hler (Proposition \ref{prop2}), by
substituting to $M$ the dense open subset where the Lee form $\theta$
does not vanish, we shall assume throughout this section that $D^g J$ nowhere vanishes, see
(\ref{integrable}).

For convenience, we choose a (local, normalized) orthonormal frame of
$\Lambda ^{+, 0} M$ of the form $\{ \phi, J \phi \}$, where $|\phi| =
\sqrt{2}$; such a frame will be called a {\it gauge}. Then, the
triple $\{ F, \phi, J \phi \}$ is a (local, normalized) orthonormal
frame of $\Lambda ^ + M$.

Recall that by Proposition \ref{prop1}  we have 
\begin{equation} \label{lambda0} W
^+ (\psi) =  - \frac{\kappa}{12}  \psi, \end{equation} 
for any section $\psi$ of $\Lambda ^{+, 0} M$, whereas 
\begin{equation} \label{lambda+} W ^+ (F) =
\frac{\kappa}{6} F. \end{equation}
With respect to the gauge $\{ \phi, J \phi \}$, the covariant
derivative $ D^g F$ is written as
\begin{equation}\label{DF}
D^g F = \alpha\otimes \phi + J\alpha\otimes J\phi,
\end{equation}
where 
\begin{equation} \alpha = \phi (J \theta); \end{equation}
equivalently,  
\begin{equation} \label{phialpha} \phi = -\frac{1}{|\theta|^2}\big(\alpha\wedge J\theta + J\alpha\wedge \theta \big); \ \  
J\phi = \frac{1}{|\theta|^2}\big(\alpha\wedge \theta - J\alpha\wedge
J\theta \big). \end{equation}
We also have
\begin{equation} \label{Dphi}
D^g \phi = - \alpha\otimes F  + \beta\otimes J\phi; \ \ 
D^g (J\phi)= -J\alpha\otimes F - \beta\otimes \phi,
\end{equation}
for some 1-form $\beta$.

From (\ref{DF}), we infer
\begin{eqnarray}\nonumber
(D^g)^2|_{\La^2M} F &=& ({\rm d}\alpha + J\alpha\wedge \beta)\otimes \phi + ({\rm d}(J\alpha) -\alpha\wedge \beta)\otimes J\phi \\ \nonumber
                    & =& -R(J\phi)\otimes \phi + R(\phi)\otimes J\phi.
\end{eqnarray}
Because of (\ref{lambda0}), this reduces to
\begin{equation}\label{ricci1}
\left\{
\begin{array}{c@{ = }c}
{\rm d}\alpha - \beta\wedge J\alpha \ & \frac{(\kappa -s)}{12}J\phi\\
{\rm d}(J\alpha) + \beta\wedge \alpha \ & -\frac{(\kappa -s)}{12}\phi. 
\end{array}
\right.
\end{equation}
Similarly, because of (\ref{lambda+}), we infer 
the following additional relation from (\ref{Dphi}):
\begin{equation}\label{ricci2}
{\rm d}\beta + \alpha \wedge J\alpha = - \frac{(s + 2\kappa)}{12}F.
\end{equation}
Notice that  1-forms $\alpha$ and $\beta$ are both {\it gauge dependent}; if
$$\phi' = (\cos\varphi ) \phi + (\sin\varphi )J\phi $$ 
they  transform to
$$\alpha'= (\cos\varphi) \alpha + (\sin\varphi) J\alpha; \ \ 
\beta'= \beta + {\rm d}\varphi.$$ 
 
We next introduce 1-forms $n_i, m_i, i=1,2$ by
\begin{equation}\label{Dtheta}
D^g \theta  = m_1\otimes \theta + n_1\otimes J\theta + m_2 \otimes \alpha + 
n_2\otimes J\alpha.
\end{equation}

By (\ref{DF}) and (\ref{phialpha}) we derive
\begin{equation}\label{DJtheta}
\begin{array}{c@{}c}
&D^g (J\theta) = -n_1\otimes \theta  + m_1\otimes J\theta -(n_2+J\alpha)\otimes\alpha +(m_2 + \alpha)\otimes J\alpha; \\ 
& \ \ \ \ \ D^g \alpha \  = -m_2\otimes \theta + (n_2 + J\alpha)\otimes J\theta +
 m_1\otimes\alpha - (n_1 -\beta)\otimes J\alpha; \\
& D^g(J\alpha) = - n_2\otimes \theta  -(m_2 + \alpha)\otimes J\theta + 
(n_1-\beta)\otimes \alpha + m_1\otimes J\alpha.
\end{array} 
\end{equation}

A  straightforward computation, using identities (\ref{ricci1}) and the fact that 
the vector field $K= ({\kappa}^{-\frac{1}{3}}J\theta)^{\sharp_g}$, the
dual  of 
${\kappa}^{-\frac{1}{3}}J\theta$, is  
Killing (see Proposition \ref{prop2}), gives the following 
expressions 
for $m_i$ and $n_i$:
\begin{equation}\label{mn}
\left\{
\begin{array}{c@{ = }c}
m_1 & m_0 + (p-\frac{(\kappa -s)}{24|\theta|^2} +\frac{1}{2}) \theta \\ 
n_1 & Jm_0  +(p-\frac{(\kappa -s)}{24|\theta|^2} -\frac{1}{2}) J\theta \\ 
m_2 & J\phi(m_0) -(p +\frac{(\kappa -s)}{24|\theta|^2} +\frac{1}{2})\alpha  \\ 
n_2 & -\phi(m_0) -(p +\frac{(\kappa -s)}{24|\theta|^2} +\frac{1}{2})J\alpha,
\end{array}
\right.
\end{equation}
where $p$ is a  smooth function, and $m_0$ is a 1-form which 
belongs to the distribution 
${\cal D}^{\perp}= {\rm span} \{ \alpha , J\alpha \}$, 
the orthogonal
complement of 
${\cal D} = {\rm span} \{ \theta , J\theta   \}$. 

Since 
 $m_1 = {\rm d}\ln |\theta|$, the 1-form $m_0$ is nothing else 
than the projection of 
${\rm d}\ln |\theta|$  to the subbundle ${\cal D}^{\perp}$. Moreover,
with respect to any gauge $\phi$,  we write
\begin{equation}\label{m0}
m_0 = q\alpha + r J\alpha,
\end{equation}
for some smooth functions $q$ and $r$. 

In view of (\ref{integrable}),  
identities (\ref{Dtheta}) and (\ref{mn}) are  
conditions on the 2-jet of $J$. Since
$J$ is completely determined by $W^+$ (see Proposition \ref{prop1}), these 
are the conditions on the 4-jet of 
the metric referred to in the introduction.

This completes the analysis of the Einstein condition and we are now going to
see how the vanishing of $W^-$ interacts on further jets of $g$.

For that,  
we introduce the ``mirror  frame'' of $\La^-M$:
$${\bar F}= -F + \frac{2}{|\theta|^2}\theta\wedge J\theta; \ \
{\bar \phi} = \phi + \frac{2}{|\theta|^2} J\alpha\wedge \theta;$$
$$I{\bar \phi}= J\phi + \frac{2}{|\theta|^2} J\alpha\wedge J\theta, $$
where the {\it negative} almost Hermitian structure $I$, of which the
anti-self-dual 2-form ${\bar F}$ is the K{\"a}hler form, is 
equal to $J$ on ${\cal D}$ and $-J$ on ${\cal D}^{\perp}$. 
By (\ref{DJtheta}) and the fact that $\theta= \frac{{\rm d}\kappa}{3\kappa}$,
we obtain the following expression for the covariant derivative of the
Killing vector field 
$K=(\kappa^{-\frac{1}{3}}J\theta)^{{\sharp_g}}$
\begin{equation}\label{DK}
D^g K = \kappa^{-\frac{1}{3}}|\theta|^2\big(q{\bar \phi} - rI{\bar\phi} - (p-\frac{1}{2}){\bar F} + \frac{(\kappa -s)}{24|\theta|^2}F\big).
\end{equation}
Moreover, since $K$ is Killing, we have
\begin{equation} \label{killing} D^g_X \Psi = R(K,X), \end{equation}
where $\Psi = D^g K$.

Considering the ASD parts of both sides  of (\ref{killing}), we infer
that the condition $W^-=0$ is equivalent to 
\begin{equation}\label{W^-=0}
D^g(\Psi^-) = \frac{s}{24}(\bar{\phi}(K)\otimes {\bar \phi} + I\bar{\phi}(K)\otimes I\bar{\phi} + IK\otimes {\bar F}),
\end{equation}
where 
$$\Psi^-= \kappa^{-\frac{1}{3}}|\theta|^2\big(q{\bar \phi} - rI{\bar\phi} - (p-\frac{1}{2}){\bar F}\big)$$ is the ASD part of $\Psi=D^gK$, 
see (\ref{DK}). 
Furthermore, by (\ref{Dtheta}) and (\ref{DJtheta}) one gets
\begin{eqnarray}\nonumber
D^g {\bar F} &=&  -(2m_2 + \alpha)\otimes {\bar \phi} + (2Jm_2 +J\alpha) \otimes I{\bar \phi};\\ 
\label{I}
D^g {\bar \phi} &=& \ \ \  (2m_2 +\alpha) \otimes {\bar F} + (2n_1- \beta)\otimes I{\bar \phi}; \\ 
\nonumber
D^g I{\bar \phi} &=& -(2Jm_2+ J\alpha)\otimes {\bar F} - (2n_1 -\beta)\otimes {\bar \phi}.
\end{eqnarray} 
Keeping in mind that 
$\theta=\frac{{\rm d}\kappa}{3\kappa}$ and $m_1={\rm d}\ln|\theta|$,
(\ref{W^-=0}) then reduces to 
\begin{eqnarray}\label{system1}
{\rm d}p  &=& -(p-\frac{1}{2})(2m_1 - \theta)  + q(m_2+ \alpha) \\ \nonumber
    & & + r(Jm_2 + J\alpha) - 
\frac{s}{24|\theta|^2} \theta \\ \label{system2} 
{\rm d}q &=& -(p-\frac{1}{2})(m_2+ \alpha) -  q(2m_1 - \theta) \\\nonumber
    & & - r(2n_1 - \beta) - \frac{s}{24|\theta|^2} \alpha \\\label{system3}
{\rm d}r &=& -(p-\frac{1}{2})(Jm_2+ J\alpha) + q(2n_1 - \beta) \\\nonumber
    & & - r(2m_1 - \theta)- 
\frac{s}{24|\theta|^2} J\alpha.
\end{eqnarray}
Now, taking into account (\ref{ricci1}) and (\ref{ricci2}), (\ref{system1})--(\ref{system3}) constitute 
a closed differential system that a self-dual Einstein Hermitian metric 
must satisfy; by  (\ref{ricci1}), (\ref{ricci2}), (\ref{DJtheta}) 
and  (\ref{mn})
one can directly check
that the integrability conditions ${\rm d(d}p)={\rm d(d}q)={\rm d(d}r)=0$ are satisfied. This 
is a first evidence that the existence of self-dual 
Einstein Hermitian metrics with prescribed 4-jet at a given point
can be expected. 
To carry out this program  explicitly, we first consider the case when
$q\equiv 0, r\equiv 0$ and show that it  precisely corresponds
to {\it cohomogeneity-one} self-dual Einstein Hermitian metrics.

\subsection{Self-dual Einstein Hermitian metrics of cohomogeneity one}

\vspace{0.4cm}

A Riemannian 4-manifold $(M,g)$ is said to be (locally) {\it of 
cohomogeneity one}, if it admits a (local)  
isometric action of a Lie group $G$, 
with three-dimensional orbits. 
The manifold $M$ is then locally a product
$$ M \cong (t_1,t_2)\times G/H.$$
The metric $g$ descends to a left invariant metric  $h(t)$ 
on each orbit $\{ t \} \times G/H$, and, by an appropriate choice of
the parameter $t$,   can be written as 
$$g= dt^2 + h(t).$$
If, moreover,  $(M,g)$ is Einstein and self-dual, and $G$ is at least
of dimension four, then, according to a result of A. Derdzi\'nski \cite{De}, 
the spectrum of the self-dual Weyl tensor 
of $g$ is everywhere degenerate, and 
$g$ is Hermitian with respect
some invariant complex structure. 

Here is a way of constructing such metrics, all belonging  to the class of {\it diagonal
Bianchi metrics of type A} (see e.g. \cite{Tod1}).  Let 
$\widetilde{G}$ be one of
the following six three-dimensional Lie groups: 
${\Bbb R}^3$,  
${\rm Nil}^3, {\rm Sol}^3$, Isom(${\Bbb R}^2$), SU(1,1) or
SU(2); let $H$ be a discrete subgroup of $\widetilde{G}$ and
consider, on $\widetilde{G} / H$, the family of diagonal metrics $h (t)$ of the form
\begin{equation}\label{diagonal}
 h(t) = A(t)\sigma_1^2 + B(t)\sigma_2^2 + C(t)\sigma_3^2, 
\end{equation}
where $A,B,C$ are positive smooth functions, 
and $\sigma_i$ are the standard left invariant
generators of the corresponding Lie algebras; we thus have
$$d\sigma^1 = n_1 \sigma_2\wedge \sigma_3; \  d\sigma_2=-n_2\sigma_1\wedge\sigma_3; \   d\sigma_3 = n_3\sigma_1\wedge \sigma_2$$
for a triple $(n_1,n_2,n_3)$, $n_i \in \{-1,0,1\}$, depending on the
chosen group, according to the following table:

\bigskip

\begin{center}
\begin{tabular}{|c|c|c|} \hline
{\rm class}  & $ n_1 \ \ \ n_2 \ \ \ n_3 $ &  ${\widetilde G}$ \\
\hline \hline
{\rm I}      & $0 \ \ \ \ 0 \ \ \ \ 0\ $       &  ${\Bbb R}^3$ \\ \hline
{\rm II}     & $0 \ \ \ \ 0 \ \ \ \ 1\ $       &  ${\rm {Nil}^3}$\\ \hline
${\rm VI}_0$ & $1 \ \ {-1} \ \ \  0\ $       &  ${\rm {Sol}^3}$ \\ \hline
${\rm VII}_0$& $1 \ \ \ \ 1 \ \ \ \ 0\ $       &  ${\rm Isom}({\Bbb R}^2)$ \\ \hline
{\rm VIII}   & $1 \ \ \ \ 1 \  {-1}\ $       &  ${\rm SU}(1,1)$ \\ \hline
{\rm IX}     & $1 \ \ \ \ 1 \ \ \ \ 1\ $       &  ${\rm SU(2)}$ \\ \hline
\end{tabular}
\end{center}
\vspace{0.4cm}

Except for Class ${\rm VI}_0$, when  $A = B$ all these metrics admit a  further (local)
symmetry which rotates the $\{\sigma_1, \sigma_2 \}$-plane, i.e. we get
the so-called {\it biaxial} Bianchi metrics, see e.g. \cite{CP}.
We thus obtain
diagonal Bianchi metrics of Class A, 
admitting a local isometric action of a four-dimensional Lee group $G$, where
$G$ is
${\Bbb R}\times{\rm Isom({\Bbb R}^2)}$, {\rm  U(1,1)}, {\rm U(2)}, or the non-trivial
central extension of ${\rm Isom}({\Bbb R}^2)$ corresponding to
biaxial Class II metrics.
Clearly, any such metric admits 
a positive {\it and}  a negative invariant
Hermitian structure, $J$ and $I$, whose K{\"a}hler forms are given by
$$F=  \sqrt{C}dt\wedge \sigma_3 + A\sigma_1\wedge \sigma_2,$$
and
$${\bar F} =\sqrt{C}dt\wedge \sigma_3 - A\sigma_1\wedge \sigma_2,$$
respectively.
When imposing the Einstein and  the self-duality conditions,   
we obtain  an ODE system 
for the unknown functions $A$ and $C$,  which can be explicitly
solved, cf. {e.g.} 
\cite{P}, \cite{Le}, \cite{DS}, \cite{Tod1}, \cite{CP},
\cite{bergery}.

In the sequel, we shall simply refer to these 
(self-dual,  Einstein,  Hermitian) metrics  as {\it diagonal Bianchi} 
metrics. 

Notice  that  $4$-dimensional locally symmetric metrics, 
  i.e. real and complex space forms, can also be put  (in several  ways) as 
diagonal Bianchi metrics. For example, self-dual Einstein Hermitian metrics in 
Class I are all flat \cite{Tod1}. 

Our next result shows that,  apart from locally symmetric
spaces, diagonal Bianchi metrics in the above sense  are actually {\it
  all}
(non-locally symmetric) cohomogeneity-one self-dual Einstein Hermitian 
metrics, and, in fact, 
can be characterized by the property 
$m_0\equiv 0$ in the notation  of the preceding section. More
  precisely, we have:

\begin{theo}\label{th2}    
Let $(M,g)$ be a self-dual
Einstein 4-manifold. Suppose that $(M,g)$ is not locally symmetric. Then the following three conditions are equivalent:
\begin{enumerate}
\item[(i)] $(M,g)$ is of  cohomogeneity one and 
the spectrum of $W^+$ is degenerate.
\item[(ii)] $(M,g)$ admits a local isometric action of a Lie group of
  dimension at least four, with three-dimensional orbits, and  is 
locally isometric to a  diagonal Bianchi self-dual Einstein Hermitian
metric belonging to one of the classes ${\rm II}$, ${\rm VII}_0$, ${\rm VIII}$ or
${\rm IX}$. 
\item[(iii)] $(M,g)$ admits a positive, non-K{\"a}hler
Hermitian structure $J$, and a negative 
Hermitian structure $I$ such that $I$ is equal to $J$ on 
${\cal D}={\rm span} \{\theta , J\theta \}$ and to $-J$ on the orthogonal 
complement ${\cal D}^{\perp}$ ; equivalently, the 1-form $m_0$
of $(g,J)$ vanishes identically.

\end{enumerate}
\end{theo}
\begin{proof}

${\rm (i)} \Rightarrow {\rm (iii)}$. By Propositions \ref{prop1} 
and \ref{prop2},  $W^+$ has two distinct, non-constant eigenvalues at
any point  and  there exists a positive, non-K{\"a}hler Hermitian structure 
$J$ whose K{\"a}hler form $F$ 
generates the eigenspace of $W^+$ corresponding to the 
simple eigenvalue. It follows that the Hermitian structure is preserved by the 
action of $G$, and therefore both functions 
$|D^g F|^2 = 2|\theta|^2$ and 
$|W^+|^2=\frac{{\kappa}^2}{24}$ are constant along the orbits of $G$; in particular, 
${\rm d}\ln|\theta|$ is colinear to $\theta = \frac{{\rm
    d}\kappa}{3\kappa}$, at any point; this means that  
$m_0=0$; by (\ref{I}) and (\ref{mn}), 
the vanishing of $m_0$ is equivalent to the integrability of the 
negative almost Hermitian structure $I$.

\vspace{0.2cm}
${\rm (iii)} \Rightarrow {\rm (ii)}$.  
If $m_0\equiv 0$  or,  equivalently,  if the negative almost Hermitian
structure $I$ is integrable, then,   by (\ref{I}), the Lie form
$\theta_I$ of $(g,I)$ reads:
\begin{equation}\label{thetaI}
\theta_I =(2p +\frac{(\kappa -s)}{12|\theta|^2})\theta.
\end{equation}
According to (\ref{mn}) we also have  
$m_1={\rm d}\ln|\theta| =(p-\frac{(\kappa -s)}{24|\theta|^2} +
\frac{1}{2}) \theta$ and  $\theta= \frac{1}{3}{\rm d}\ln|\kappa|$; 
it follows that ${\rm d}\theta_I=0$;  then, locally, 
$\theta_I = {\rm d}f$ for a positive function $f$, i.e., $g$ is 
conformal to a K{\"a}hler metric ${g'} = f^2g$.
Since $W^-=0$, the K{\"a}hler metric $g'$ is 
of zero scalar curvature. Clearly, the Killing
field $K$ preserves both $J$ and $g$, hence, also, 
the K{\"a}hler structure $(g',I)$. Two cases occur, according as $g '$ is
homothetic or not to $g$.

\vspace{0.2cm}
(a) Suppose ${g'}$ is {\it not} homothetic to $g$; equivalently,  the scalar 
curvature $s$ of $g$ does not vanishes; then, by \cite{De}, 
$K'=I{\rm grad}_g(f^{-1})$ is a Killing vector field for $g$ and $g'$ 
and is 
holomorphic with respect $I$. By the very definition of $I$ we have
that 
$J|_{\cal D} = I|_{\cal D}$; the Killing
vector fields $K'$ and $K$ are thus colinear everywhere (see
(\ref{thetaI}));  
it follows that
$K'$ is a 
constant multiple of $K$. By considering  $z=f^2$ as a local
coordinate on $M$ and,  by introducing a holomorphic
coordinate $x+iy$ on the (locally defined) 
orbit-space for the holomorphic action of $K+\sqrt{-1}IK$ on
$(M,I)$,  
the metric $g$ can be  written in the following form:
\begin{equation}\label{g}
g = \frac{1}{z^2}[e^{u}w({\rm d}x^2 + {\rm d}y^2) + w {\rm d}z^2 + w^{-1}\om^2],
\end{equation}
where $u(x,y,z)$ is a smooth function satisfying the 
${\rm SU(\infty )}$ Toda field equation: 
$$u_{xx} + u_{yy} + (e^u)_{zz}=0,$$
$w$ is a positive function given by 
$$w= \frac{6(zu_z -2)}{s},$$
and $\om$ is a 
connection 1-form of the  ${\Bbb R}$-bundle 
$M \mapsto N=\{(x,y,z)\} \subset {\Bbb R}^3$, 
whose curvature is given by
\begin{equation}\label{domega}
{\rm d}\om = -w_x{\rm d}y\wedge {\rm d}z - w_y {\rm d}z\wedge {\rm d}x - (we^u)_z 
{\rm d}x\wedge {\rm d}y, 
\end{equation}
(see, e.g.   \cite{Tod2}). 
Moreover, the Killing field $K$ is dual to $\frac{1}{wz^2}\om$, and
the (anti-self-dual) K{\"a}hler form of 
the negative Hermitian structure $I$  is given by
\begin{equation}\label{barF}
{\bar F} =\frac{1}{z^2}\big(we^u {\rm d}x\wedge {\rm d}y - {\rm d}z\wedge \omega \big).
\end{equation}
By (\ref{thetaI}) we have that 
 ${\cal D}={\rm span}\{ \theta , J\theta \}= 
{\rm span}\{\theta_I, I\theta_I \}={\rm span}\{K^{\sharp_g},IK^{\sharp_g}\}$, so that
the K{\"a}hler form $F$ of the positive Hermitian structure $J$ is given by
\begin{equation}\label{F}
F = \frac{1}{z^2}\big(we^u {\rm d}x\wedge {\rm d}y + {\rm d}z\wedge \omega \big).
\end{equation}
It is now easily seen that (\ref{barF}) and (\ref{F}) 
simultaneously define integrable almost complex structures if and only if 
$w_x=w_y=0$, or equivalently if and only if
$u(x,y,z)=u_1(x,y) + u_2(z)$. This means that $u$ is a {\it separable} solution to 
the ${\rm SU(\infty )}$ Toda field equation. Up to a  change of the
holomorphic coordinate $x+ iy$, it 
is explicitly given by \cite{Tod2}
$$e^u = \frac{4(c + bz + az^2)}{(1 + a(x^2 + y^2))^2},$$
for properly chosen constants $a,b,c$.
Any such  solution gives rise to a
{diagonal Bianchi} self-dual Einstein Hermitian metric pertaining to 
one of classes II, ${\rm VII}_0$, VIII and IX,
depending on the choice of the constants
$a,b,c$ (see e.g. \cite[Sec. 8]{CP}) for a common case of these
metrics in the Bianchi IX case).

\vspace{0.2cm}
(b) If $g ' $ is homothetic to $g$, i.e. $(g,I)$ is itself a K{\"a}hler structure of zero scalar curvature, 
then $g$ is locally hyperk{\"a}hler  and $K$
is a Killing vector field preserving the K{\"a}hler structure $I$. Then,
one of the two 
following situations occurs:
 
(b1) {\it $K$ is {\it triholomorphic}}, i.e.  $K$
preserves each  K{\"a}hler structure in the hyperk{\"a}hler family: Then
the quotient space,  $N$,  for the (real) action of $K$ is flat and is
endowed with a field of parallel straight lines. 
This situation is described by the Gibbons-Hawking Ansatz \cite{GH},
and the metric $g$ has  the form: 
$$g = w({\rm d}x^2 + {\rm d}y^2 + {\rm d}z^2) + \frac{1}{w}\omega^2,$$
for a positive harmonic function $w(x,y,z)$ on $N$ and a 1-form $\omega$ 
on $M$ satisfying
$${\rm d}\om = -w_x{\rm d}y\wedge {\rm d}z - w_y {\rm d}z\wedge {\rm d}x - w_z 
{\rm d}x\wedge {\rm d}y.$$
The Killing field $K$ is dual to $\frac{1}{w}\om$ 
and one may consider
that
the positive and negative Hermitian structures, $J$ and $I$, 
correspond to the 2-forms
$$F= w{\rm d}x\wedge {\rm d}y + {\rm d}z\wedge \omega; \ \
{\bar F}= w{\rm d}x\wedge {\rm d}y - {\rm d}z\wedge \omega,$$ 
respectively. We again conclude  $w_x=0, w_y=0$, and therefore
$w=az +b$. The case $a=0$ corresponds to flat metrics in Class I,
whereas,  when
$a\neq 0$,  by putting 
$at = az + b, \sigma_1 = {\rm d}x, \sigma_2 = {\rm d}y, \sigma_3 = \omega$,  
the metric becomes a  diagonal Bianchi metric of Class II.

(b2) {\it $K$ is {\it not} triholomorphic}: 
Since, nevertheless,  $K$ preserves $(g,I)$, 
the metric $g$ takes the form \cite{BoF}
$$ g = e^{u}w({\rm d}x^2 + {\rm d}y^2) + w {\rm d}z^2 + w^{-1}\om^2, $$
where $u(x,y,z)$ is a solution to the 
${\rm SU(\infty)}$ Toda field equation,  $w=au_z$, 
$\omega$ satisfies (\ref{domega}) and  
$a$ is a constant. Moreover, $K$ is dual to $\frac{1}{w}\omega$, and $I$ is defined by the anti-self-dual form 
$${\bar F} = we^u {\rm d}x\wedge {\rm d}y - {\rm d}z\wedge \omega.$$
Similar arguments as above show that 
$w_x=w_y=0$, i.e., $u$ is  a separable solution to  
the ${\rm SU(\infty )}$ Toda field equation, and therefore
our metric is again a diagonal Bianchi 
metric in one of the classes II, ${\rm VII}_0$, VIII or IX, cf.
\cite{CP}.

\vspace{0.2cm}
The implication ${\rm (ii)} \Rightarrow {\rm (i)}$ is clear.
\end{proof}

\begin{rem}{\rm
A  weaker version of Theorem \ref{th2} was announced in \cite{de} 
(see \cite[Rem.~1.3]{de} and Lemma 2 above).}
\end{rem}

\subsection{The generic case} We now consider the generic case, when $m_0$  a {\it non-vanishing} section of 
${\cal D}^\perp$, hence determines a gauge $\phi$ such that  
$r\equiv 0, q\neq 0$ in (\ref{mn}). According to
(\ref{mn}), the 1-form $\alpha$ is then given by
\begin{equation}\label{alpha}
m_1 = {\rm d}\ln|\theta| = q\alpha + (p- \frac{(\kappa -s)}{24|\theta|^2} + \frac{1}{2})\theta;
\end{equation}
moreover,   by (\ref{system1})--(\ref{system3}), we have that 
\begin{eqnarray}\label{beta}
\beta &= & \frac{1}{q}\Big(p(2p + \frac{(\kappa -s)}{12|\theta|^2} - 1) - \frac{\kappa}{24|\theta|^2} + 2q^2\Big)J\alpha \\\nonumber
      &  &-\frac{(\kappa -s)}{12|\theta|^2}J\theta,
\end{eqnarray}
\begin{eqnarray}\label{frobenius1}
{\rm d}p &=&\Big(2q^2 - p(2p - \frac{(\kappa -s)}{12|\theta|^2} - 1) 
       - \frac{\kappa}{24|\theta|^2}\Big)\theta \\ \nonumber
   & & - q\Big(4p +\frac{(\kappa -s)}{12|\theta|^2} - 1\Big)\alpha,
                   \\\label{frobenius2} 
{\rm d}q &=& - q\Big(4p -\frac{(\kappa -s)}{12|\theta|^2} - 1\Big)\theta 
         \\ \nonumber 
   & & - \Big(2q^2 -p(2p + \frac{(\kappa -s)}{12|\theta|^2} - 1) 
  + \frac{\kappa}{24|\theta|^2}\Big)\alpha.
\end{eqnarray}
By differentiating (\ref{alpha}) and  by making use of
(\ref{frobenius1})--(\ref{frobenius2}), 
we get
\begin{equation}\label{dalpha}
{\rm d}\alpha = \frac{(\kappa -s)}{12|\theta|^2}\alpha \wedge \theta =
\alpha\wedge J\beta;
\end{equation}
this is nothing else than  the first relation in (\ref{ricci1}),  when 
$\beta$ is given by (\ref{beta});
by substituting the expression (\ref{beta}) for $\beta$ 
into the second relation of (\ref{ricci1}),  we obtain
\begin{equation}\label{dJalpha}
{\rm d}(J\alpha) = J\alpha \wedge J\beta.
\end{equation}
In view of (\ref{alpha}) and (\ref{frobenius1})--(\ref{frobenius2}),
it is not hard to check  that  the 1-form $J\beta$ is equivalently given by
\begin{equation}\label{dJbeta}
J\beta = {\rm d}\ln(\frac{|\kappa|}{|q||\theta|^4}),
\end{equation}
so that  (\ref{dJalpha}) becomes
\begin{equation}\label{dJa}
{\rm d}(\frac{\kappa}{q|\theta|^4}J\alpha)=0;
\end{equation}
from (\ref{DJtheta}) we get
\begin{equation}\label{dJthe}
{\rm d}(J\theta)= J\theta\wedge\big(\frac{1}{3}{\rm d}\ln|\kappa| -2{\rm d}\ln|\theta|\big) + J\alpha \wedge \eta,
\end{equation}
or, equivalently, 
\begin{equation}\label{dJtheta}
{\rm d}(\frac{\kappa^{\frac{1}{3}}}{|\theta|^2} J\theta)= 
\frac{\kappa^{\frac{1}{3}}}{|\theta|^2}J\alpha\wedge \eta,
\end{equation}
where
$$\eta = -2q\theta + (2p + \frac{(\kappa -s)}{12|\theta|^2} -1)\alpha.$$

We are now rea{\rm d}y to prove  the existence of self-dual Einstein Hermitian metrics 
with $m_0\neq 0$. More precisely, we  exhibit a 1--1-correspondence
between these metrics and the set of solutions of the integrable Frobenius system 
(\ref{frobenius1})--(\ref{frobenius2}).
We start with the data 
$(s, \kappa, |\theta|)$  consisting 
of a constant $s$ (the scalar curvature), a 
nowhere vanishing smooth function $\kappa$ (the conformal scalar curvature),
and a positive smooth function $|\theta|$ (the norm of the Lie form
$\theta= \frac{{\rm d}\kappa}{3\kappa}$), defined on an open subset ${U}$ 
of $M$, 
such that $\theta\wedge {\rm d}|\theta|^2$ has no zero on ${U}$ 
(equivalently, $m_0$ does not vanish on ${U}$).
We then introduce local
coordinates $x= \kappa^{\frac{1}{3}} \neq 0$ and $y=|\theta|^2 >0$. 
Observe that $x$ is a {\it momentum map} for the Killing field $K$ 
with respect to the self-dual K{\"a}hler metric 
${\bar g}={\kappa}^{\frac{2}{3}}g$  while $y=|K|_{\bar g}^2$ is the 
square-norm of $K$ with respect to ${\bar g}$ (see Proposition \ref{prop2}). 
The Lee form $\theta$ is then given by
\begin{equation}\label{deftheta} 
\theta= \frac{{\rm d}x}{x},
\end{equation}   
and the 1-form $\alpha$ is 
given by (\ref{alpha}) for some smooth functions $p(x,y)$ and 
$q(x,y)\neq 0$ of $x,y$, i.e. 
\begin{equation}\label{defalpha}
\alpha = \frac{1}{q}\Big( \frac{{\rm d}y}{2y} - \frac{1}{x}(p -\frac{(x^3-s)}{24y} + \frac{1}{2}){\rm d}x\Big).
\end{equation}
Then, (\ref{frobenius1})--(\ref{frobenius2}) can be made into the
following  
Frobenius system for the 
(unknown) functions $p$ and $q^2$:
\begin{eqnarray}\label{defp}
{\rm d}p &=& \frac{1}{x}\Big[ 2q^2 + 2(p+\frac{(x^3 -s)}{24 y})(p-\frac{(x^3 -s)}{24 y} +1) -\frac{1}{2} - \frac{x^3}{24y}\Big]{\rm d}x \\ \nonumber
   & & - \frac{1}{y}\Big[2p + \frac{(x^3 -s)}{24y} - \frac{1}{2}\Big]{\rm d}y
\end{eqnarray}
\begin{eqnarray}\label{defq}
{\rm d}(q^2) &=& -\frac{1}{y}\Big[ 2q^2 -2p(p+\frac{(x^3 -s)}{24y} - \frac{1}{2}) +\frac{x^3}{24y}\Big]{\rm d}y\\ \nonumber
& &- \frac{2}{x}\Big[\Big(p-\frac{(x^3 -s)}{24y} + \frac{1}{2}\Big)\Big(2p(p +\frac{(x^3-s)}{24y} -\frac{1}{2}) - \frac{x^3}{24y}\Big)\\ \nonumber
 & & \ \ \ \ \ \ \ - 2q^2(1-p)\Big]{\rm d}x
\end{eqnarray}
A straightforward computation shows 
that the integrability condition ${\rm d}({\rm d}p)={\rm d}({\rm
  d}q^2)=0$ is satisfied (as a matter of fact, the explicit solutions are given in Lemma 3
below). The above mentioned 
correspondence
between solutions to (\ref{defp})--(\ref{defq}) and 
self-dual Einstein Hermitian metrics with $m_0\neq 0$ now goes as
follows. 
Since  (\ref{defp})--({\ref{defq}) is integrable, 
each  value  of $(p,q)$ at a given point $(x_0,y_0)$ can be extended
to a solution of  (\ref{defp})--(\ref{defq})
in some neighborhood 
$V$ of 
$(x_0,y_0)$; moreover, by choosing  $q(x_0,y_0) \neq 0$,  we may assume
that $q$ has no zero  
on $V$; 
by (\ref{defalpha}) and (\ref{defp})--(\ref{defq}),  one 
immediately obtains (\ref{dalpha}) for the corresponding 1-form
$\alpha$. We then 
introduce a third local coordinate, $z$, such that
\begin{equation}\label{defJalpha} 
J\alpha = 
\frac{qy^2}{x^3}{\rm d}z,
\end{equation} 
see (\ref{dJa}). Finally, since the 1-form $J\theta$ 
satisfies 
(\ref{dJthe}) or,  equivalently, (\ref{dJtheta}), 
the integrability condition reads as follows:
$${\rm d}(\frac{qy}{x^2}\eta)=0,$$  
see (\ref{dJa}) and (\ref{dJthe}); by using
(\ref{frobenius1})--(\ref{dJalpha}),  one easily checks that 
the integrability condition  is actually satisfied, so that 
\begin{equation}\label{defJtheta}
J\theta= \frac{y}{x}({\rm d}t + h{\rm d}z),
\end{equation} where $t$ is a suitable 
transversal coordinate to $(x,y,z)$,
and $h(x,y)$ is a smooth function on $V$, defined by
$${\rm d}h =-\frac{qy}{x^2}\eta.$$
It is an easy consequence of (\ref{defp}) that the above equation 
is solved  by
\begin{equation}\label{defh}
h = \frac{yp}{x^2} + \frac{x}{24}. 
\end{equation}
The metric $g$ and the orthogonal
almost complex structure $J$ are then given  by
$$g = \frac{1}{|\theta|^2}(\theta\otimes\theta + J\theta\otimes
J\theta + \alpha\otimes \alpha + J\alpha\otimes J\alpha);$$
according to (\ref{deftheta}),(\ref{defalpha}),(\ref{defJalpha}) and
(\ref{defJtheta}), and by using the  coordinates 
$(x,y,z,t)$,  the metric $g$ takes the form 
\begin{equation}\label{canonic}
g = \frac{1}{y}\Big[\frac{{\rm d}x^2}{x^2} + \frac{1}{q^2}\Big(\frac{{\rm d}y}{2y} - \frac{1}{x}(p -\frac{(x^3 -s)}{24y} + \frac{1}{2}){\rm d}x\Big)^2 + \frac{q^2y^4}{x^6}{\rm d}z^2 +
\frac{y^2}{x^2}({\rm d}t + h{\rm d}z)^2\Big]; 
\end{equation}
this shows that any self-dual Einstein Hermitian metric 
with $m_0\neq 0$ is 
locally isometric to a metric of  the above form  for some solution $(p,q)$ to 
(\ref{defp})--(\ref{defq}). 

Conversely, for any solution to
(\ref{defp})--(\ref{defq}), the corresponding almost-Hermitian metric 
$(g,J)$ is 
self-dual Einstein Hermitian metric with $m_0\neq 0$.  Indeed,
by (\ref{dalpha}), (\ref{dJalpha}) and (\ref{dJtheta}), 
$J$ is integrable and it is easily checked that $\theta=\frac{{\rm d}x}{x}$ is the Lee form for $(g,J)$, 
i.e., 
$${\rm d}F = -2\theta\wedge F;$$ 
moreover, the 1-form $\alpha$ corresponds to
the gauge
$$ \phi = -\frac{1}{y}\big(\alpha\wedge J\theta + J\alpha\wedge \theta \big),$$
meaning that $\alpha = \phi(J\theta)$; one directly computes 
$${\rm d}\phi = (\theta + J\beta)\wedge \phi,$$
where the 1-form $\beta$ is given by (\ref{beta}); it follows that
$\beta$ is precisely the {\rm 1-form} defined by (\ref{Dphi}) and that 
(\ref{dalpha})--(\ref{dJalpha}) are nothing else than the 
Ricci identities (\ref{ricci1});  this allows us  to recognize the curvature:
By (\ref{ricci1}),  the Ricci tensor of $(g,J)$ is 
$J$-invariant, and, since $\theta= \frac{{\rm d}x}{x}$, 
the dual vector field $K$ of 
$\kappa^{-\frac{1}{3}}J\theta=\frac{1}{x}J\theta$  
is Killing, cf. e.g. \cite{AG}; by (\ref{dJtheta}) and (\ref{DF}),  
the covariant derivative of 
$\theta$ is  given by (\ref{Dtheta}) for $p$ and $q$ constructed as above,
and $r\equiv 0$; hence, (\ref{beta}) and (\ref{frobenius1})--(\ref{frobenius2})
(equivalently, (\ref{defp})--(\ref{defq}))
are  the same as relations
(\ref{system1})--(\ref{system3}); these, in turn, are  a way of
re-writing (\ref{W^-=0}); it follows that
the projection of the curvature to $\La^-M$ reduces to
$\frac{s}{12}{\rm Id}|_{\Lambda^-M}$, i.e. the Hermitian metric $g$ 
is Einstein and self-dual,  with scalar curvature equal to $s$, 
see (\ref{SO(4)}); turning back to (\ref{dalpha}), 
we conclude that the conformal scalar curvature is $\kappa=x^3$, 
see (\ref{ricci1}); the metric constructed in this way is not of
cohomogeneity one,  as 
$m_0\neq 0$, see 
Theorem \ref{th2}. Finally, 
different solutions $(p,q)$ of  (\ref{defp})--(\ref{defq}) give rise
to non-isometric metrics,  
as $p$ and $q$ are completely determined by $|W^+|, {\rm d}|W^+|$ and ${\rm d}|D^gW^+|$, 
see Sec.~ 2 and (\ref{alpha}). 

We finally observe that the metric (\ref{canonic}) admits two commuting vector fields, 
$\frac{\partial}{\partial t}$ and $\frac{\partial}{\partial z}$.

We summarize the results obtained so far as  follows:
\begin{theo}\label{th3} Let $(M,g,J)$ be a self-dual Einstein Hermitian 
4-manifold. Suppose that  $(M,g,J)$  is neither locally-symmetric  nor
of  cohomogeneity one. 
Then,  on an open dense subset of $M$,  $g$ is locally given by (\ref{canonic}). In particular, $(M,g)$ admits  a 
local isometric action
of  ${\mathbb R}^2$ almost-everywhere.
\end{theo}

\begin{rem} \label{rem3} {\rm (i) It is easily seen that the metrics (\ref{canonic})
have only 
2-dimensional continuous symmetries. Moreover, as we already observed,  
the coordinate $x=\kappa^{\frac{1}{3}}$
is  a {momentum map}  of the Killing vector field 
$\frac{\partial}{\partial t}$ with respect to the 
K{\"a}hler metric 
${\bar g}= x^2 g$ while, by (\ref{defp}) and (\ref{defh}),  
a momentum map $\tilde{\mu}$ of the second Killing field, 
$\frac{\partial}{\partial z}$, 
is given by
$$2x{\tilde \mu} = y + \frac{x^3 + s}{12},$$
where $\frac{x^3 + s}{12} = 
\frac{\kappa + s}{12}$ is the (pointwise constant)
holomorphic sectional curvature of $(g,J)$. 

The momentum map $x$ is also equal to the scalar curvature
of the K{\"a}hler metric ${\bar g}$. 
A straighforward computation shows that the second momentum map $\tilde{\mu}$
defined above is related to the Pfaffian  of the {\it normalized Ricci form}
${\bar \sigma}$ of the K{\"a}hler metric ${\bar g}$ by 
$$ \tilde{\mu} = 12 \, ( {\rm Pfaff} \, {\bar \sigma}  + b),$$
where $b$ is the constant appearing in (\ref{explicitef}) below. This
fits with an observation of R. Bryant in \cite{Br}.
({\rm Recall that for any
  $2$-form $\psi$, the Pfaffian of $\psi$ with respect to ${\bar g}$
  is defined by: $\psi \wedge \psi = 2\,  {\rm Pfaff}\,  \psi \, v _{{\bar
      g}}$, where $v _{{\bar
      g}}$ is the volume form of ${\bar g}$; the
  normalized Ricci form ${\bar \sigma}$ is the $(1,1)$-form associated
  to the normalized Ricci tensor,  ${\bar S}$,  appearing in the
  usual decomposition ${\bar R} = {\bar S} \wedge {\bar g} + W$ of the
  curvature operator of ${\bar g}$~; it is related to the usual Ricci
  form ${\bar \rho}$ by ${\bar \sigma}  = \frac{1}{2} \, ({\bar \rho} _0 +
  \frac{x}{12} \, {\bar \omega})$, where ${\bar \rho} _0$ is the
  trace-free part of ${\bar \rho}$; since $g = x ^{-2} {\bar g}$ is
    Einstein and ${\rm d}^c x$ is the dual of a Killing vector field,
    we have that  ${\bar \rho} _0 = - \frac{1}{x} \, ({\rm d}
  {\rm d}^c x) _0$; the result follows easily}).

(ii) It follows from Theorems \ref{th2} and \ref{th3} that 
every self-dual Einstein 
Hermitian  4-manifold admits a (local) isometric ${\Bbb R}^2$-action 
compatible 
with a {\it product structure} in the sense of \cite{joyce}; 
the general considerations in 
\cite[Sec.2]{joyce} therefore apply to the present  situation;  
a  detailed analysis
of self-dual Einstein 4-manifolds admitting ${\Bbb R}^2$-continuous symmetry 
has been carried out by D. Calderbank \cite{Ca0}, based on results of \cite{Ca1}. }
\end{rem}

\vspace{0.2cm}
We end this  section by providing an explicit form for the metric
(\ref{canonic}), in view of the following
\begin{Lemma}\label{integrate}  
The solutions $p(x,y)$ and $q(x,y)$ of the system (\ref{defp})--(\ref{defq})  are explicitly given by
\begin{equation}\label{explicitep}
p = \frac{f}{y^2} - \frac{(x^3 -s)}{24y} + \frac{1}{4};
\end{equation}
\begin{equation}\label{expliciteq}
q^2= \frac{1}{y^2}\Big[\frac{x}{2}f' -f + \big(\frac{x^3-s}{24}\Big)^2\Big] - \frac{x^3}{24y} - p^2,
\end{equation}
where  
\begin{equation}\label{explicitef}
f(x)= ax^2 + bx^4 - \frac{(x^6 - s^2)}{576},
\end{equation}
 $a$ and $b$ are constants defined by positivity in (\ref{expliciteq}), 
and $f'$ stands for the first derivative of $f$.
\end{Lemma}
\begin{proof}
We first observe that (\ref{defp}) can be equivalently written as
$${\rm d}\Big(y^2(p + \frac{(x^3 -s)}{24y} - \frac{1}{4})\Big) = $$
$$\frac{y^2}{x}\Big[2q^2 + 2(p + \frac{(x^3 -s)}{24y})(p - 
\frac{(x^3 -s)}{24y}) + 2(p + \frac{(x^3 -s)}{24y} - \frac{1}{4}) + \frac{x^3}{12y}\Big]{\rm d}x;$$
this shows that $y^2(p + \frac{(x^3 -s)}{24y} - \frac{1}{4})$ is
function of $x$, say $f$; from the above equality,  we get (\ref{explicitep}) and 
(\ref{expliciteq}), where $f$ is a (still unknown)  smooth function;
in order to determine $f$, we differentiate (\ref{expliciteq}) by
using (\ref{explicitep}) and  
substitute into (\ref{defq}); then, cancellations occur and  (\ref{defq}) 
eventually reduces to 
\begin{equation}\label{ODE}
x^2f'' -5xf' + 8f + \frac{(x^6 -s^2)}{72}=0; 
\end{equation}
the solutions of (\ref{ODE}) are given 
by (\ref{explicitef}). \end{proof}

\section{Self-dual Einstein Hermitian metrics with hyperhermitian structures}

In this section, we consider self-dual,  Einstein,  Hermitian metrics
which in addition admit a {\it non-closed} hyperhermitian structure compatible 
with the negative orientation. It is well-known that 
LeBrun-Pedersen metrics, which are of  cohomogeneity one under the action of the
unitary group ${\rm U(2)}$, carry such hyperhermitian structures; 
in LeBrun's coordinates \cite{Le} these metrics read as follows:
\begin{equation}\label{g1}
g = \frac{1}{(b t^2 + 4c)^2}
\Big((1 + \frac{8b}{t^2} + \frac{16c}{t^4})^{-1}{\rm d}t^2 + \frac{t^2}{4}
\big[\sigma_1^2 + \sigma_2^2 + (1 + \frac{8b}{t^2} + \frac{16c}{t^4})
\sigma_3^2 \ \big]\Big),
\end{equation}
where $b$ and $c$ are properly chosen constants \cite{Mad}; more
precisely, we have the following 
\begin{prop}\label{prop5}{\rm (\cite{Mad})}
Let $(M,g)$ be an oriented self-dual 
Einstein 4-manifold. Assume that $(M,g)$ admits  a  $\Lie{U}(2)$ isometric action with
generically three-dimensional $\Lie{SU}(2)$-orbits. 
If $g$ admits a non-closed, $\Lie{U}(2)$-invariant  
negative hyperhermitian structure, then 
$g$ is isometric to (\ref{g1}) with $c > b^2$,
and  actually admits exactly two distinct invariant hyperhermitian 
structures. 
\end{prop} 
We here prove the following more general result:
\begin{theo}\label{th1} A self-dual Einstein Hermitian 4-manifold
$(M,g,J)$ locally admits a non-closed, negative hyperhermitian structure 
if and only if
$g$ is  locally isometric to one of the $\Lie{U}(2)$-invariant metrics
(\ref{g1}) with $c > b^2$; then, $(M, g)$  actually carries 
exactly two distinct hyperhermitian structures, each of them $\Lie{U}(2)$-invariant.
\end{theo}
 
We first establish general facts concerning  self-dual Einstein 4-manifolds 
which carry a {\it non-closed} hyperhermitian structure compatible with the negative orientation. 
As already observed in Sec.2, a (negative) 
hyperhermitian structure $(g, I_1,I_2,I_3)$ is determined by a real $1$-form  
$\theta$ --- the common Lee form of $(g,I_i)$, also the Lee form of
the Obata connection --- satisfying conditions 
(\ref{EW}) and (\ref{hypherm}), and such that $\Phi := {\rm d} \theta$
is self-dual; in particular, the 2-form $\Phi$ is harmonic. 
The next Lemma 
shows that the self-dual Weyl tensor of $g$ is completely determined by 
$\theta$, $\Phi$ and the 
first covariant derivative $D^g\Phi$ of $\Phi$.
\begin{Lemma}\label{gauduchon1} Let $(M,g)$ be  an oriented self-dual 
Einstein 4-manifold and assume that $(M, g)$ carries a 
negative hyperhermitian structure. Then, as a symmetric operator
acting on $\Lambda ^ + M$,  the self-dual Weyl tensor $W^+$ is 
given by
\begin{equation}\label{gau1}
W^+(\psi) = \frac{1}{2}[\psi, \Phi] + \frac{1}{|\theta|^2} D^g_{\psi(\theta)} \Phi,
\end{equation}
where $\psi$ is any self-dual 2-form, $\theta$ is  viewed 
as a vector field by Riemannian duality,
and $[\cdot, \cdot]$ denotes the commutator of 
2-forms,  viewed as skew-symmetric endomorphisms of the tangent bundle.
Moreover, $\theta$ and $\Phi$ are related by
\begin{equation}\label{gau3}
D^g_{\theta} \Phi = 2|\theta|^2\Phi.
\end{equation}
\begin{equation}\label{gau2}
{\rm d}|\theta|^2 -(\frac{s}{12} + |\theta|^2)\theta + \Phi(\theta)=0,
\end{equation}
\end{Lemma}
\begin{proof} By using (\ref{EW}), the right-hand side of 
$$R_{X,Y}\theta= (D^g)^2_{Y,X}\theta - (D^g)^2_{X,Y}\theta $$
is easily computed; we thus obtain:
\begin{eqnarray}\label{util1}
R(\theta\wedge Z) &=& -\frac{1}{2}{\rm d}|\theta|^2\wedge Z - \frac{1}{2}(\frac{s}{12} -|\theta|^2)\theta \wedge Z\\ \nonumber
& &  
-\frac{1}{2}\Phi(Z)\wedge \theta - \frac{1}{2}D^g_Z\Phi + \theta(Z)\Phi.
\end{eqnarray}
Since $g$ is self-dual and Einstein, $R = \frac{s}{12}{\rm Id}|_{\Lambda^2M} + W^+$, 
see (\ref{SO(4)}). 
Then, by projecting (\ref{util1}) to $\La^-M$, we get  
(\ref{gau2}), whereas  
the projection of (\ref{util1}) to $\La^+M$ gives (\ref{gau1}) and (\ref{gau3}). 
\end{proof}
\begin{cor}\label{cor1} {\rm (\cite{ET,Ca})} Every hyperhermitian 
structure on a conformally flat 
4-manifold is closed.
\end{cor}
\begin{proof}   If we assume  that 
$\Phi \neq 0$ somewhere on $M$ and that the anti-self-dual Weyl tensor 
is identically zero, then,  
after contracting (\ref{gau1}) and (\ref{gau3}) with $\Phi$, we obtain 
$\theta = \frac{1}{4} {\rm d}\ln|\Phi|^2$, 
which contradicts $\Phi = {\rm d}\theta \neq 0$. \end{proof}

We can compute the covariant derivative $D^g_{\theta}W^+$ of $W^+$ 
along the dual vector field of $\theta$ (still denoted by $\theta$), by using (\ref{gau1})   together with 
(\ref{gau3}) and (\ref{gau2}) (the latter are used for evaluating the term 
$(D^g)^2_{\theta,\psi(\theta)} \Phi$ which  appears in the
calculation); we thus get  
\begin{Lemma}\label{gauduchon2} Let $(M,g)$ be an oriented self-dual 
Einstein 4-manifold, admitting a 
negative hyperhermitian structure; then,   the covariant derivative $D^g_{\theta} W^+$ of the 
self-dual Weyl tensor $W^+$ along the dual vector field of the
Lee form $\theta$ is given by
\begin{eqnarray}\nonumber  
\big( (D^g_{\theta} W^+)(\psi), \phi \big) &=& \big( [W^+(\phi),\psi] + [W^+(\psi),\phi], \Phi \big) \\ \nonumber
& & + (4|\theta|^2 - \frac{s}{6})\big( W^+(\psi), \phi \big)   \\ \label{gau4}
& & + |\Phi|^2\big( \psi, \phi \big) - 
3\big( \Phi, \psi \big) \big( \Phi, \phi \big) ,
\end{eqnarray}
for any sections, $\phi$ and $\psi$, of $\La^+M$.
\end{Lemma}

From Lemma \ref{gauduchon2} and Propositions \ref{prop1} and
\ref{prop2},  we infer
\begin{prop}\label{prop6} Let $(M,g)$ be an oriented self-dual 
Einstein 4-manifold, admitting a non-closed
hyperhermitian structure compatible with the negative orientation. 
Then the following three conditions are equivalent:
\begin{enumerate}
\item[{\rm(i)}] the spectrum of $W^+$ is everywhere degenerate;
\item[{\rm(ii)}] $W^+$ has  two distinct eigenvalues at any point; 
\item[{\rm(iii)}] the self-dual 2-form $\Phi$ is a nowhere vanishing
eigenform for $W^+$ with respect to the simple eigenvalue,  
and is proportional to a positive Hermitian structure $J$.
\end{enumerate}\
\end{prop}
\begin{proof}

${\rm (i)} \Rightarrow {\rm (ii)}$. 
According to Proposition \ref{prop1}, if  the spectrum of $W^+$ is
everywhere degenerate, then  either $W^+$  vanishes identically (and therefore
 the hyperhermitian
structure is closed by Corollary \ref{cor1})  or
$W^+$ has two distinct eigenvalues $\la$ and $-\frac{\la}{2}$ at any
point.

${\rm (ii)} \Rightarrow {\rm (iii)}$.
By Proposition \ref{prop1}, we know  that
a normalized generator $F$ of the $\la$-eigenspace of $W^+$ is 
the K{\"a}hler form of a positive Hermitian structure $J$. 
Let $\phi$ be any self-dual 2-form orthogonal to 
$F$, with $|\phi | ^2 = 2$; then,  $\phi$ and $\psi = (J\circ \phi)$ are orthogonal, 
$(-\frac{\la}{2})$-eigenforms of $W^+$; by substituting into
(\ref{gau4}),  we get
$$ 0= \big( (D^g_{\theta} W^+)(\phi), \psi \big) = - 3\big( \Phi, \psi \big) \big( \Phi, \phi \big),$$
$$ -{\rm d}\la (\theta)= \big( (D^g_{\theta} W^+)(\phi), \phi \big) = -(4|\theta|^2 -\frac{s}{6})\la +2|\Phi|^2 - 3\big( \Phi, \phi \big)^2,$$
$$ -{\rm d}\la (\theta)= \big( (D^g_{\theta} W^+)(\psi), \psi \big) = -(4|\theta|^2 -\frac{s}{6})\la + 2|\Phi|^2 - 3\big( \Phi, \psi \big)^2.$$
From the last two equalities,  we get
$\big( \Phi, \psi \big) = \pm \big( \Phi, \phi \big) $, and 
by the first one we  
conclude that $\big( \Phi, \psi \big) = \big( \Phi, \phi \big) =0$. 
This shows that $\Phi$ is a multiple of $F$. It remains to prove that
$\Phi$ does not vanish on $M$; by taking a two-fold cover of $M$ if
necessary, we may assume that the Hermitian structure $J$ is globally
defined on $M$; by 
 Proposition \ref{prop2}, $(g,J)$ is conformally K{\"a}hler
and $\la^{\frac{2}{3}}F$ 
is the corresponding 
closed K{\"a}hler form; but $\Phi$ is also a closed, self-dual 2-form, 
and a multiple of ${F}$, hence a 
constant (non zero) multiple of 
$\la^{\frac{2}{3}}F$.

${\rm (iii)} \Rightarrow {\rm (i)}$. 
This is an immediate consequence of Proposition \ref{prop1}.
\end{proof}

\vspace{0.2cm}
\noindent
{\bf Convention:}
From now on,  we assume that $(M,g)$ is an oriented self-dual
Einstein 4-manifold 
whose self-dual Weyl $W^+$ has degenerate spectrum, and  which 
admits a {\it non-closed} hyperhermitian structure compatible with the 
negative orientation of $M$. According to Proposition \ref{prop6},
$W^+$ has two distinct eigenvalues which we denote by $\la$ and $-\frac{\la}{2}$, 
and  the harmonic self-dual 2-form $\Phi$ 
defines a positive Hermitian structure $J$ on 
$(M,g)$  whose K{\"a}hler form, $F$, is an $\la$-eigenform for $W^+$.  Moreover,
it follows from Proposition \ref{prop2} that,  
after rescaling the metric if necessary, we may assume:
\begin{equation}\label{util2}
\Phi = \frac{1}{2}\la^{\frac{2}{3}}F.
\end{equation}
In the notation  of Sec.2.1, the conformal 
scalar curvature $\kappa$ of $(g,J)$ is thus equal to $6\la$; 
the Lee form
$\theta_J$ and the Killing vector field $K$, rescaled by an appropriate 
positive constant,  are therefore given by:
\begin{equation}\label{theta-X}
\theta_J= \frac{{\rm d}\la}{3\la}; \ \ K=J\rm{grad_g}(\la^{-\frac{1}{3}}),
\end{equation} (see Proposition \ref{prop2}).

At this point, our main technical result reads as follows: 
\begin{prop}\label{prop7}  A self-dual 
Einstein Hermitian 4-manifold $(M,g,J)$  admits a 
non-closed, hyperhermitian structure compatible with the negative orientation 
if and only if
the Lee form $\theta_J$ satisfies
\begin{equation} \label{gau5} \begin{split} 
D^g \theta_J &= \frac{(1+ \la^{\frac{2}{3}})(s +
  3\la^{\frac{1}{3}})}{12}g \\ & \ \ \ + 
\frac{(1 + 2\la^{\frac{2}{3}})}{(1 + \la^{\frac{2}{3}})} \theta_J\otimes \theta_J + 
\frac{\la^{\frac{2}{3}}}{(1+\la^{\frac{2}{3}})}J\theta_J\otimes J\theta_J.
\end{split} \end{equation}
In this case, 
$(M,g)$ actually admits exactly two non-closed hyperhermitian structures 
$\{ I_1',I_2',I_3' \}$ and $\{ I_1'', I_2'', I_3'' \}$ whose Lee forms, 
$\theta '$ and $\theta ''$, are given  by 
$$\theta ' = \frac{1}{(1 + \la^{\frac{2}{3}})}\big( \theta_J -  \la^{\frac{1}{3}}
J\theta_J \big),$$
$$\theta '' = \frac{1}{(1 + \la^{\frac{2}{3}})}\big( \theta_J + \la^{\frac{1}{3}}
J\theta_J \big)$$
respectively.
Moreover, the Killing vector field $K$ is triholomorphic 
for both hyperhermitian structures, i.e.,
$K$ preserves all   complex structures  $I_i'$ and $I_i''$, $i=1,2,3$.
\end{prop}
\begin{proof}

We first  show that  
if $(M,g,J)$ admits a non-closed hyperhermitian structure compatible with 
the negative 
orientation, then the corresponding  Lee form $\theta$ must be one of the forms $\theta'$ and $\theta''$ 
given in Proposition \ref{prop7}. 

From (\ref{gau3}) and the fact
that $\Phi$ is an $\la$-eigenform of $W^+$,  we infer 
\begin{equation}\label{util3}
{\rm d}|\Phi|^2 =4|\Phi|^2\theta + 4\la \Phi(\theta).
\end{equation}
By differentiating (\ref{util3}) and  by using (\ref{gau2}) in order to
compute  
${\rm d}(\Phi(\theta ))$,  we obtain
$$({\rm d}\la - 3\la \theta)\wedge \Phi(\theta ) + \big(|\Phi|^2 + \la(\frac{s}{12} + 
|\theta|^2)\big)\Phi =0; $$
we infer:
\begin{equation}\label{util4}
|\Phi|^2 = - \la(\frac{s}{12}+ |\theta|^2).
\end{equation}
By substituting the above expression of  $|\Phi|^2$ in (\ref{util3}),
and by using  (\ref{gau2}) again, we get
\begin{equation}\label{util5}
{\rm d}\la - 3\la\theta = \frac{3\la^2}{|\Phi|^2}\Phi(\theta).
\end{equation}
Now, according to the above convention,  
by (\ref{theta-X}) and (\ref{util2}) we end up with the following
expression for $\theta$:
\begin{equation}\label{util7}
\theta = \frac{1}{(1 + \la^{\frac{2}{3}})}\big( \theta_J -  \la^{\frac{1}{3}}
J\theta_J \big).
\end{equation}
This shows that every non-closed  hyperhermitian structure 
is completely determined by the self-dual harmonic 2-form $\Phi$. 
It remains to prove that 
$\Phi$ itself is 
determined, up to sign,  by the metric $g$; then, the two possible values of 
$\theta$ appearing in Proposition \ref{prop7} will only differ  by conjugation of 
$J$ or,   equivalently, by substituting $ - \Phi$ to  $\Phi$.
Notice  that, according to  our convention, at this stage we have
the freedom to rescal the $2$-form $\Phi$ by a non-zero constant.   
In other words, by 
fixing one non-closed hyperhermitian structure and by following our 
convention, we know that  
any other non-closed hyperhermitian structure corresponds 
to a harmonic 2-form of the form
$a \Phi = \frac{a}{2}\la^{\frac{2}{3}}F$, where $a$ is a non-zero
constant. Our claim is that $a=\pm 1$; to see this,
by using  (\ref{gau1}) and (\ref{gau3}),  we calculate
$$|D^g \Phi |^2 = 2|\theta|^2(3|\Phi|^2 + |W^+|^2);$$
in the present situation, when $W^+$ has degenerate spectrum,   
the norm of $W^+$ is given by $|W^+|^2=\frac{3}{2}\la^2$;  
then, by (\ref{util4}), 
the above equality reduces itself to
\begin{equation}\label{util6}
|D^g \Phi |^2 = -(\frac{|\Phi|^2}{\la} + \frac{s}{12})(6|\Phi|^2 + 3\la^2);
\end{equation} 
it is readily checked that if the $2$-forms
$\Phi$ and $a \Phi$ simultaneously satisfy (\ref{util6}), then $a=\pm 1$. 

\vspace{0.2cm}
We now check that the conditions (\ref{EW})\&(\ref{hypherm}) for either 
$\theta'$ or $\theta''$ 
are equivalent 
to (\ref{gau5}).
Keeping   (\ref{util2}) in mind, we see that (\ref{util5}) can be equivalently 
re-written as
\begin{equation}\label{util8}
\theta_J= \theta + \la^{\frac{1}{3}}J\theta;
\end{equation}
then, the equivalence 
``(\ref{gau5}) $\Leftrightarrow$ (\ref{EW})\&(\ref{hypherm})'' follows by 
a straightforward computation involving the expressions  
(\ref{util7}) and (\ref{util8}), and using formula (\ref{integrable});
the 1-forms 
$\theta'$ and $\theta''$ thus 
correspond to two distinct, non-closed hyperhermitian structures 
$\{ I_1',I_2',I_3' \}$ and $\{ I_1'', I_2'', I_3'' \}$ provided that 
(\ref{gau5}) holds, see Sec. 1.2.

\vspace{0.2cm} 
As a final step, we have to prove that $K$ is triholomorphic with
respect to both
hyperhermitian structures. For a general hyperhermitian 
structure $I_i, i=1,2,3$,  with Lee form $\theta$, and 
for any  Killing field $K$,  we have
$${\cal L}_K I_i = D_K I_i -[DK,I_i],$$
where $D$ is the
Weyl derivative given by (\ref{D^J});  we thus only need  to check 
that in our specific situation $D K$ commutes with $I_i$; by using (\ref{D^J}), (\ref{theta-X}),
(\ref{integrable}) and
(\ref{gau5}),   we get 
$$DK = \theta(K){\id} + \frac{(1+ \la^{\frac{2}{3}})}{4} J;$$
the claim follows immediately. \end{proof}

\begin{cor}\label{cor2} {\rm (\cite{ET})} A locally-symmetric self-dual 
Einstein 4-manifold does not admit non-closed hyperhermitian structures.
\end{cor}
\begin{proof} Any such manifold 
is either a space of constant curvature, hence conformally flat, or 
a K{\"a}hler manifold of constant holomorphic sectional curvature (see 
Propositions \ref{prop1} and \ref{prop2}). In the former  case,  the claim follows 
by Corollary \ref{cor1}, whereas  in the latter case $\theta_J=0$;  we
then 
conclude by using   
Proposition \ref{prop7}. \end{proof}

\begin{rem} {\rm D. Calderbank  proved that
any conformal selfdual 4-manifold admitting two distinct Einstein-Weyl structures
is equipped with a canonical conformal submersion to an Einstein-Weyl 
3-manifold \cite{Ca2}. 
In the situation described by Proposition \ref{prop7},  this conformal
submersion is seen as follows: the  hyperhermitian structures $\{
I_1',I_2',I_3' \}$ and $\{ I_1'', I_2'', I_3'' \}$ determine a SO(3)-valued
function, $p$, on $M$ defined by: 
$$I_i'' = \sum_{j=1}^3 a_{ij}I_j'; \ A=(a_{ij})\in {\rm SO(3)};$$
we claim that $p$ is a conformal submersion of $(M,g)$ to 
SO(3)=${\Bbb RP}^4$: The
differential of $p$ is easily computed by using the fact that $I_i''$ and $I_j'$ are both  integrable; we thus obtain: 
\begin{equation}\label{dA}
{\rm d}(a_{ij}) + \frac{\la^{\frac{2}{3}}}{2(1+ \la^{\frac{2}{3}})}\Sigma_{k=1}^3
a_{ik}\big([I'_k,I'_j] K\big)^{\sharp_g}=0;
\end{equation}
here,  $[\cdot, \cdot]$ denotes the commutator of endomorphisms of $TM$
and $^{{\sharp}_g}$ stands for the Riemannian duality; from
(\ref{dA}), we infer:  
$${\cal L}_{K} a_{ij} =0,$$
$$\sum_{i,j} \big(da_{ij}(X)\big)^2 =
\frac{\la^{\frac{4}{3}}}{2(1+\la^{\frac{2}{3}})^2}g(X,X), \ \forall 
X\in K^{\perp};$$
The first equality shows 
that $p$ coincides with the  projection of $M$ to
the space,  $N$, of orbits of $K$, whereas 
the second equality means  that the $K$-invariant metric 
${\bar g}=\frac{\la^{\frac{2}{3}}}{(1+\la^{\frac{2}{3}})}g$ descends
to the round  
metric of  ${\rm SO(3)} = {\mathbb R} P^3$; in other words,  $K$ defines a Riemannian submersion from $(M,{\bar g})$
to ${\rm SO(3)}$.}
\end{rem}

\vspace{0.2cm}
\noindent
{\bf Proof of Theorem \ref{th1}.} 
We first notice that 
the Killing vector field $K$ is trivial 
if and only if $\la$ is constant (see (\ref{theta-X})), or, equivalently, 
$\theta_J=0$. 
Thus, according to Propositions \ref{prop6} and \ref{prop7}, 
if $(M,g,J)$ is a self-dual Einstein Hermitian 4-manifold 
admitting a {\it non-closed} hyperhermitian structure,  
the Killing vector field 
$K$ does not vanish
on an open, dense subset of $M$. 
It then follows from \cite{GT,CT,CP} that self-dual Einstein 4-manifolds
admitting two distinct hyperhermitian structures and a non-trivial 
triholomorphic  Killing 
vector field are locally  given by 
Proposition \ref{prop5}.

For completeness, however, we here give a different and 
more direct argument adapted to our ``Hermitian'' situation.

By Proposition \ref{prop5} it is sufficient to show that our metric can be 
written
in  the diagonal form (\ref{diagonal}). Since the eigenvalues of $W^+$ are 
not constant, i.e., $\theta_J\neq 0$ (Proposition \ref{prop7}), 
we introduce the variable $t=\la^{\frac{1}{3}}$; 
the Lee form $\theta_J$ is then equal to  $\frac{{\rm d}t}{t}$, whereas  
the dual $1$-form of the Killing vector field is given by $-\frac{1}{t^2}J{\rm d}t$. 
We set: $\sigma_3 = f(t)J{\rm d}t$,
for some smooth function $f$ of $t$, 
and we insist that 
\begin{equation}\label{dsigma3}
{\rm d}\sigma_3 = \sigma_1\wedge \sigma_2,
\end{equation} 
where the 1-forms $\sigma_1$ and $\sigma_2=J\sigma_1$ are 
both orthogonal to ${\rm d}t$  and satisfy
\begin{equation}\label{dsigma1}
{\rm d}\sigma_1 =\sigma_2\wedge \sigma_3; \ \ {\rm d}\sigma_2 = \sigma_3\wedge \sigma_1.
\end{equation}
We then derive  $f$ from (\ref{dsigma3}): 
By differentiating (\ref{util7}) and by making use of (\ref{util2}),   we obtain
\begin{equation}\label{dJr}
{\rm d}(J{\rm d}t) = -\frac{(1+t^2)t^2}{2}F + \frac{2t}{(1+t^2)}{\rm d}t\wedge J{\rm d}t.
\end{equation}
By   (\ref{util8}), (\ref{util4}) and (\ref{util2}),  we also get
$$|{\rm d}t|^2 =-(\frac{t}{2} + \frac{s}{12})(t^4 + t^2);$$
it follows that $\big({\rm d}\sigma_3, {\rm d}t\wedge J{\rm d}t \big) =0$ if and only 
if $(\ln f)'= -\frac{2t}{(1+t^2)} - \frac{1}{(t + \frac{s}{6})}$, where the 
prime 
stands for $\frac{{\rm d}}{{\rm d}t}$; we then have   
$f= \frac{a}{(1+t^2)(t+ \frac{s}{6})}$, hence 
\begin{equation}\label{sigma3}
\sigma_3 = \frac{a}{(1+t^2)(t + \frac{s}{6})} J{\rm d}t
\end{equation} 
for a positive constant $a$.

In order to determine the 1-forms $\sigma_1$ and $\sigma_2$,  we choose 
a gauge 
$\phi$ or,  equivalently, a  1-form $\alpha =\phi(J\theta_J) \in {\cal
  D}^{\perp}$; since $\sigma_1$ and $\sigma_2=J\sigma_1$ are orthogonal to ${\rm d}t$, 
there certainly 
exists  
a smooth function $h$ of $t$ and a smooth function $\varphi$ on $M$, such that
$$\sigma_1 = h(\cos\varphi \alpha + \sin\varphi J\alpha); 
\sigma_2=h(-\sin\varphi \alpha + \cos\varphi J\alpha);$$
by (\ref{sigma3}) and (\ref{dsigma3}),  we obtain the following
expression for $h$:
\begin{equation}\label{h}
h^2 = \frac{at^2}{(t + \frac{s}{6})^2(1+t^2)};
\end{equation} 
by using (\ref{sigma3}) and (\ref{ricci1}), 
we now see that the conditions (\ref{dsigma1}) 
are equivalent to
\begin{equation}\label{varphi}
{\rm d}\varphi + \beta  + \frac{(\frac{s}{6}-t^3 +at)}{t(1+t^2)(\frac{s}{6}+ t)}J{\rm d}t =0;
\end{equation}
therefore, the existence of a smooth function $\varphi$ satisfying
(\ref{varphi}) is equivalent to the following condition:
$${\rm d}(\beta  + \frac{(\frac{s}{6}-t^3 +at)}{t(1+t^2)(\frac{s}{6}+ t)}J{\rm d}t)=0;$$
a straightforward computation involving (\ref{ricci2}) and (\ref{dJr}) shows
that the above equality holds whenever the constant $a$ is chosen
equal to $1 + \frac{s^2}{36}.$  \ \

\section{Hermitian structures on quaternionic quotients}

Let $(N,g)$ be a quaternionic K{\"a}hler manifold of real dimension $4n$, 
endowed with a non-trivial Killing field $K$ which preserves  the
quaternionic structure. According to 
Galicki \cite{galicki1, galicki2} and
Galicki-Lawson \cite{G-L}, under some  ``non-degeneracy'' condition
for $K$ one can define
a $4(n-1)$-dimensional quaternionic orbifold $(M,g^*)$
via the so-called {\it quaternionic reduction construction}. This
can be described as follows. 
We  first consider the following orthogonal
splitting 
of the bundle of 2-forms: 
\begin{equation}\label{quaternion-split}
\La^2N = \La^+N \oplus \La^{1,1}N \oplus \La^{\perp}N,
\end{equation}
where: 
\begin{enumerate}
\item[$\bullet$] $\La^+ N$ is the  3-dimensional 
sub-bundle of ``self-dual'' 
2-forms which determines the {quaternionic structure} (also identified
to 
a sub-bundle $A^+N$ of skew-symmetric endomorphism
of $TN$):  both $A^+N$ and $\La^+N$ are preserved by
the Levi-Civita connection,  $D ^g$, and at each point $x$ of $N$
there is an orthonormal basis $\{ I_1, I_2, I_3 \}$ of $A^+N \subset {\rm End}(T_x N)$ with
the property that:  $I_i\circ I_j = -\delta_{ij}{\rm Id}|_{TN} + 
\epsilon_{ijk} I_k$ 
(resp. $\La^+N = {\rm span}(\omega_1,\omega_2,\omega_3)$, 
where $\omega_i$ are the fundamental 2-forms  of the almost Hermitian 
structures $(g,I_l)$. In the sequel, we  refer to any such choice of 
$I_l$'s (resp. $\omega_l$'s) as a {\it trivialization} of $A^+ N$
(resp. $\La^+ N$);
\item[$\bullet$] $\La^{1,1}N$  is the sub-bundle of 2-forms which are 
$I_i$-invariant for any section of $A^+N$;
\item[$\bullet$] $\La^{\perp}N$ denotes the orthogonal complement of 
$\La^+N \oplus \La^{1,1}N$ in $\La^2N$.
\end{enumerate}
We denote by ${\Pi^+}$ the projection of $\La^2N$ to $\La^+N$; 
for any trivialization $\{ \omega_1, \omega_2, \omega_3 \}$ of 
$\La^+ N$ we then have
$$\Pi^+ = \frac{1}{2n}\sum_{l} \omega_l\otimes \omega_l,$$ and 
$\Pi^+_{K} :=
\frac{1}{2n}\sum_{l} (i_K\omega_l \otimes \omega_l)$ is a section of
$T^*N\otimes \La^+ N$. Then, Galicki-Lawson showed  \cite[Th.~2.4]{G-L}.
that there exists
a section $f_K$ of $\La^+N$ such that
$${\rm d}^{D^g}  f_K = D ^g  f_K = \Pi^+_K.$$
The section $f_K$ is called {\it the momentum map} associated to
$(N,g,K)$ and it is easily seen that the ``level set''
\begin{equation}\label{LK}\nonumber
L_{K} := \{ x\in N: f_K(x)=0 \}
\end{equation}
is $K$-invariant. 

Assuming that $K_x \neq 0$ at $x\in L_{K}$, 
Galicki-Lawson proved that $L_{K}$ is regular, 
i.e. $L_K$ is a smooth submanifold of $N$. If  moreover
the quotient space $M:= L_{K}/K$ is (locally) a $(4n-4)$-dimensional
manifold (or just an orbifold), then it becomes a 
quaternionic K{\"a}hler manifold with respect to the ``projected'' 
quaternionic structure, $g^*$,  of $N$.
Thus, when $N$ is 8-dimensional, the quaternonic reduction gives rise to a 
four dimensional {\it anti-self-dual} Einstein orbifold
(with respect to the canonical orientation induced by $N$). 
Note that when $K$ is the generator of a
$S^1$-quaternionic action on $N$, under
the non-degeneracy condition as above $M$ always inherits an orbifold
structure, cf. \cite[Th. 3.1 \& Cor. 3.2]{G-L}.

The above construction applies in particular to 
$N = {\mathbb H}{P}^2$ endowed with certain {\it weighted} $S^1$-actions;
one thus obtains a wealth of examples of
{\it compact} anti-self-dual Einstein orbifolds; as shown by
Galicki-Lawson, the corresponding orbifolds are all 
weighted projective planes ${\mathbb C} P^{[p_1,p_2,p_3]}$ for
some integers $0<p_1\le p_2\le p_3$ satisfying $p_3<p_1+p_2$,
\cite[Sec. 4]{G-L}.
Notice  that,  with respect
to the orientation induced by the canonical complex structure, the  metric 
becomes {\it self-dual}. (In the case when $p_1=p_2=p_3$ one obtains
the Fubini-Study metric on ${\mathbb C}{P}^2$). 

On the other hand, R. Bryant showed \cite[Sec. 4.2]{Br}
that each weighted projective plane admits a
self-dual K{\"a}hler metric which under the above assumption
for the weights
has everywhere positive scalar curvature. Therefore, according to
\cite[Lemma \ref{de-ga}]{AG},
Bryant's metric
gives rise to
a self-dual {\it Einstein} Hermitian metric on
${\mathbb C} P^{[p_1,p_2,p_3]}$, $p_3<p_1+p_2$. 

When considering  both results together, a natural question arises:

\vspace{0.2cm}
\noindent
{\bf Question.} \cite{LeBrun}
Are the Galicki-Lawson metrics on ${\mathbb C} P^{[p_1,p_2,p_3]}$
Hermitian
with respect to some anti-self-dual complex structure?

\vspace{0.2cm}

In this section we show that this is indeed the case, at least on a
dense open subset; more generally,
we show that 
the answer to the above question is essentially yes  for any anti-self-dual
Einstein 4-orbifold obtained
by quaternionic reduction from  the 8-dimensional Wolf spaces
${\mathbb H}P^2$, $SU(4)/S(U(2)U(2))$ and the corresponding
non-compact dual spaces (but according to \cite{kris}
the  argument fails for quaternionic quotients of the exeptional 8-spaces
$G_2/SO(4)$ and $G^2_2/SO(4)$). More precisely, we have the following

\begin{prop}\label{quat-quot} Let $(N,g)$ be
${\mathbb H}{P}^2, SU(4)/S(U(2)U(2))$, 
or one of the corresponding non-compact dual spaces.
Then, any anti-self-dual, Einstein
4-orbifold $(M,g^*)$ which is obtained as a quaternionic reduction
of $(N,g)$ by a quaternionic Killing field $K$ 
locally admits (a negatively oriented) 
Hermitian structure $J$. In particular, the metric $g^*$ is locally given 
by the explicit constructions in Sec. 2.
\end{prop}

The proof is based on the following simple observation.
\begin{Lemma} \label{Phi} 
Let $(N,g)$ be a quaternionic K{\"a}hler manifold 
of non-zero scalar curvature and $K$ be a Killing
field on $N$. Denote by $\Psi(X,Y)=(D ^g_X K, Y)$ the
2-form corresponding to $D ^g K$ and let $\Psi^+ = \Pi^+(\Psi)$
be the projection of $\Psi$ to
$\La^+ N$.  Then, up to multiplication by a constant,
the momentum map $f_K$ of $K$  is given by $\Psi^+$.
\end{Lemma}
\begin{proof} Since $K$ is Killing, equality (\ref{killing})
$$D ^g_X \Psi = R(K\wedge X)$$
holds.
For a quaternionic K{\"a}hler manifold the curvature operator $R$
acts on $\La^+ N$ by $\lambda {\rm Id}|_{\La^+N}$,
where $\lambda$ is a positive multiple of the scalar curvature,
cf. e.g. \cite{salamon}. Thus, projecting (\ref{killing}) to
$\La^+N$ we get $D ^g_X \Psi^+ = \lambda \Pi^+_K. $
\end{proof}   

By Lemma \ref{Phi}
the ``level set'' $L_K$ of $K$ is the same as the set of points
$x\in N$ where $\Psi^+_x =0$. Thus, at any point $x\in L_K$
the tangent space $T_xL_K$ is given by
$T_xL_K = \{ T_xN \ni X : D ^g_X \Psi^+ =0 \}.$
Since by assumption $K$ does not vanish on $L_K$, we conclude 
by (\ref{killing}) and the fact that 
$R|_{\La^+N}= \lambda {\rm Id}|_{\La^+N}$ 
$$T_xL_K = {\rm span}(I_1K,I_2K,I_3K)^{\perp},$$
where $\{I_1, I_2,I_3\}$ is any trivialization of $A^+N$. 

We also observe  that the 2-form $\Psi$ is a section
of $\La^+N \oplus \La^{1,1}N$, provided that 
$K$ preserves the quaternionic structure. Indeed,
$$[D ^g K, I_l] = D ^g _K I_l  -{\cal L}_K I_l, $$
where $[\cdot ,\cdot]$ stands for the commutator of ${\rm End}(TN)$. 
Since $K$ is quaternionic, the left-hand-side of the above equality 
is a section of $\La^+ N$. By summing over $l$ in the above relation
we  get
\begin{equation}\label{la11}
\Psi + 2\Pi^{1,1}(\Psi) \in \La^+N,
\end{equation}
where 
$\Pi^{1,1}$ denotes the projection to $\La^{1,1}N$: 
\begin{equation}\label{pi11}
\Pi^{1,1}(\psi)(\cdot,\cdot) = \frac{1}{4} \Big[(\psi(\cdot, \cdot) + \sum_l\psi(I_l\cdot ,I_l\cdot)\Big], \ \forall \psi \in \La^2N.
\end{equation}
Thus, $\Psi$ is a section of
$\La^+N \oplus \La^{1,1}N$, and at $x\in L_K$, $\Psi_x$
actually belongs to  $\La_x^{1,1}N$. 

Since $\Psi = \frac{1}{2} {\rm d} K^{\sharp}$,
where $K^{\sharp}$ is the $g$-dual 1-form of $K$, we conclude that
$${\cal L}_K \Psi = {\rm d}(i_K(\Psi)) = -\frac{1}{2} {\rm d}({\rm d}|K|^2)=0, $$
i.e. $\Psi$ is a closed $K$-invariant 2-form.
This shows that $\Psi$ projects to $M= L_K/K$ to define an
{\it anti-self-dual} form  on $(M,g^*)$, then denoted by $\Psi^*$.
Considering the Riemannian submersion 
$$\pi: L_K \longmapsto M = L_K/K,$$ 
the {\it horizontal} space, $H$, of $TL_K$ is given by
$$H = {\rm span}(K,I_1K,I_2K,I_3K)^{\perp}.$$ 
Note that $H$ is $I_l$-invariant for any section $I_l$ of $A^+N$.
Using the above remarks we calculate:
\begin{equation}\label{important}
(D ^{g^*}_{U^*} \Psi^*)(V^*, T^*) = (D ^g_U \Psi)(V,T)
 -\frac{4}{|K|^2_g}\Pi^{1,1}(i_U\Psi \wedge i_K\Psi)(V,T),
\end{equation}
where $D ^{g^*}$ is the Levi-Civita connection of $g^*$, $U^*,V^*,T^*$ 
are any vectors on $M$, and $U,V,T$ are the corresponding horizontal lifts.

By assumption,  $K$ has no zero on $L_K$; 
it then follows from (\ref{important}) and  (\ref{killing}) that
$\Psi^*$ does not  vanish identically on $M$. Thus,
on the open subset of $(M,g^*)$ where $\Psi^* \neq 0$ the
normalised ASD form $\frac{{\sqrt 2} \Psi^*}{|\Psi^*|_{g^*}}$
determines  a
{\it negative} almost Hermitian structure $J$. 
By virtue of the Riemannian Goldberg-Sachs (\cite[Prop. 1]{AG}), 
Proposition \ref{quat-quot} follows from the following
\begin{Lemma}\label{integrab} The almost-complex structure $J$ is integrable.
\end{Lemma}
\begin{proof} 
We denote  $Z^*_i$ any complex (1,0)-vector field of $(M,J)$ and  $Z_i$ 
the corresponding horizontal lift (considered as complex vector in
$T_x^{\mathbb C} N$);
then, $J$ is integrable if and only if the following identity holds:
\begin{equation}\label{integrability0}
D ^{g ^*}_{Z^*_i} (\frac{{\sqrt 2}\Psi^*}{|\Psi^*|_{g^*}})(Z^*_j,Z^*_k) =
(D ^{g ^*}_{Z^*_i} \Psi^*)(Z^*_j, Z^*_k) =0 \ \forall i,j,k ;
\end{equation}
by the very definition of $J$ we have
$\Psi(Z_i,Z_j)=0$; moreover, since $\Psi$ belongs to $\La^{1,1}N$ on $L_K$,
the almost complex structure $J$ (defined on $H$)
commutes with $I_l$'s for any trivialization $\{I_1,I_2,I_3\}$ of $A^+N$.
Then, by (\ref{important}) and (\ref{killing}) it is easily seen  
that
the integrability condition (\ref{integrability0}) for $J$  is the
same as 
\begin{equation}\label{integrability}
(D ^{g ^*}_{Z^*_i} \Psi^*)(Z^*_j, Z^*_k)=
(\na_{Z_i} \Psi)(Z_j, Z_k) = (R(K\wedge Z_i),Z_j\wedge Z_k) = 0.
\end{equation}

We now derive 
(\ref{integrability}) from
the structure of the curvature tensor of the Riemannian
symmetric spaces ${\mathbb H}{P}^2, SU(4)/S(U(2)U(2))$ 
and the corresponding  non-compact duals, ${\mathbb H}{H}^2$ and
$SU(2,2)/S(U(2)U(2))$ (we refer to
\cite{salamon, gauduchon} 
for a general description of the curvature operator, $R$, of a 
Riemannian symmetric space).  

We first consider the simplest  case
of $N={\mathbb H}{P}^2 = Sp(3)/(Sp(1)Sp(2))$ (or its non-compact
dual).   The eigenspaces of $R$ are  then 
the simple 
factors ${\bf sp}(1)$ and ${\bf sp}(2)$ 
of the isotropy Lie sub-algebra ${\bf h} = {\bf sp}(1) \oplus {\bf sp}(2)$, 
and  
the
orthogonal complement ${\bf h}^{\perp}$ of ${\bf h}$ in the space 
${\rm Skew}({\bf m})$ of 
the skew-symmetric endomorphisms of ${\bf m}= {\bf sp}(3)/{\bf h}$ 
(note that $R$ acts trivially on ${\bf h}^{\perp}$); the
decomposition ${\rm Skew}({\bf m}) = {\bf sp}(1) \oplus {\bf sp}(2)
\oplus {\bf h}^{\perp}$ into eigenspaces of
$R$ then fits with the splitting (\ref{quaternion-split});
$\La^+N$ is thus identified to    
${\bf sp}(1)$, and $\La^{1,1}N$  to ${\bf sp}(2)$, whereas
$\La^{\perp}N$ 
corresponds to the kernel of $R$, the space ${\bf h}^{\perp}$. This
shows that the curvature operator
acts on the first 
two 
factors in (\ref{quaternion-split}) by multiplication with a 
non-zero constant (a certain multiple of the scalar curvature), and
acts trivially on the
third factor (therefore, $R$ has thus
three distinct eigenvalues, $\la,\mu$ and $0$); this observation also 
shows that any
Killing field on ${\mathbb H}P^2$ is necessarily quaternionic. 

As already observed, the almost complex structure $J$ 
(defined on $H$)
commutes with the $I_l$'s, so that $I_l(Z_k)$ is again a
(1,0)-vector of $(H,J)$; we thus get
$$\Pi^+(Z_j\wedge Z_k) = \sum_{l} (Z_j, I_l(Z_k))\omega_l= 0,$$
which means that $Z_j\wedge Z_k$ is an element of
$\La_x^{1,1}M \oplus \La_x^{\perp}N$. It then follows that 
\begin{eqnarray}\nonumber
(R(K\wedge Z_i), Z_j\wedge Z_k) &=& (R(Z_j\wedge Z_k), K\wedge Z_i) \\ \nonumber
   &=&
\mu(\Pi^{1,1}(Z_j\wedge Z_k), K\wedge Z_i). 
\end{eqnarray}
But $\Pi^{1,1}(Z_j\wedge Z_k)$ is again a (2,0)-vector of $(M,J)$ 
(see formula (\ref{pi11})),
so that
$(\Pi^{1,1}(Z_j\wedge Z_k), K\wedge Z_i)=0$; this implies
(\ref{integrability}).

The same argument holds for the non-compact dual space ${\mathbb H}H^2$.

The case of $N=SU(4)/S(U(2)U(2))$ (or its non-compact dual)
is similar, but $N$ is now a {\it Hermitian symmetric} space,  whose 
canonical Hermitian structure $I$ comutes with any $I_i \in \La_x^+ N$.
The corresponding
K{\"a}hler form, $\Om_I$, then belongs to the space $\La^{1,1}N$ and gives rise
to a further splitting 
$$\La^{1,1}N = {\mathbb R}\cdot \Om_I \oplus \La^{1,1}_0 N, $$
where $\La^{1,1}_0 N$ is the orthogonal complement of $\Om_I$.
Correspondingly, the eigenspaces of the curvature $R$ are the bundles
$\La^+ N$, ${\mathbb R}\cdot \Om_I$, $\La^{1,1}_0 N$, and $\La^{\perp} N$.
Note that $R$ acts trivially
on $\La^{\perp}N$, whereas $\Om_I$ is an eigenform 
of $R$ corresponding to the simple eigenvalue; in particular,
$K$ must preserve $I$ and $\Omega_I$, so that $\Psi$ is of type $(1,1)$
with respect to $I$; in other words, the almost complex structure $I$
commutes with $J$,
when acting on $H$. It  follows that $Z_i\wedge Z_j$ belongs to
$\La^{1,1}_0 N \oplus \La^{\perp}N$,
and we conclude as in the case of ${\mathbb H}P^2$.
\end{proof}

\vspace{0.2cm}
\noindent
{\bf Remark 5.} (i) By (\ref{important}) and  Lemma \ref{integrab},
we see that $\frac{1}{|K|^2} \Psi^*$ is a  harmonic 2-form
on $(M,g^*)$; it is actually the K{\"a}hler form of a self-dual 
K{\"a}hler metric  in the conformal class of $g^*$ 
(see \cite[Prop. 2]{AG}). In particular, if $(M,g^*)$ is not a 
real space form, then $\Psi^*$ has no zero on $M$. By construction,
$\frac{2}{|K|^2}\Psi^*$ is  the curvature form 
of the submersion $\pi: L_K \longmapsto M$. It follows that 
$L_K$ is a Sasakian manifold fibered over a K{\"a}hler self-dual
--- equivalently, a Bochner-flat ---
four-manifold. It is well known that the corresponding CR-structure
of  $L_K$
has vanishing 
fourth-order Chern-Moser curvature; therefore $L_{K}$ is uniformized 
over $S^5$ with respect to ${\rm Aut}_{CR}(S^5)=PU(3,1)$,  
cf. \cite{webster}.

(ii) As observed in \cite[p. 20]{G-L}, the 
quaternionic reduction procedure can be applied 
to  the quaternionic hyperbolic space to obtain
{\it smooth}, {\it complete} (non locally symmetric) 
Einstein self-dual metrics 
of negative scalar curvature, which are necessarily 
Hermitian by Lemma \ref{integrab}; 
see also \cite{Br} for another construction of complete 
Einstein self-dual Hermitian metrics. In view of our first remark,
these examples seem to contradict 
some results in \cite{kamishima}.

\end{document}